\documentclass[12pt]{amsart}
\usepackage[T1]{fontenc}
\usepackage[english]{babel}
\usepackage{color}
\usepackage{mathtools}
\usepackage{graphicx}
\usepackage{amssymb}
\usepackage{tikz}
\usetikzlibrary{arrows}
\usepackage{hyperref}
\usepackage{dsfont,soul}
\usepackage{enumitem}
\usepackage{geometry}
\newcommand{\diagh}{.7}
\setlist[enumerate]{labelsep=*, leftmargin=1.5pc,
topsep=1ex plus0.5ex minus0.2ex,
itemsep=1ex plus0.5ex minus0.2ex,
font=\rmfamily,
font=\upshape}
\setlist[itemize]{labelsep=*, leftmargin=1.5pc,
topsep=1ex plus0.5ex minus0.2ex,
itemsep=1ex plus0.5ex minus0.2ex,
font=\rmfamily,
font=\upshape}
\newtheorem{thm}{Theorem}[section]
\newtheorem{cor}[thm]{Corollary}
\newtheorem{lem}[thm]{Lemma}

\newtheorem*{thm*}{Theorem}
\theoremstyle{definition}
\newtheorem{defn}[thm]{Definition}
\newtheorem{fact}[thm]{Fact}
\newtheorem{exa}[thm]{Example}
\newtheorem{rem}[thm]{Remark}
\numberwithin{equation}{section}
\newcommand{\bA}{\mathbb A}
\newcommand{\bC}{\mathbb C}
\newcommand{\bE}{\mathbb E}

\newcommand{\bN}{\mathbb N}

\newcommand{\bR}{\mathbb R}

\newcommand{\cC}{\mathcal C} % dual variety -> eliminate
\newcommand{\cE}{\mathcal E}
\newcommand{\cF}{\mathcal F}

\newcommand{\cL}{\mathcal L}
\newcommand{\cM}{\mathcal M}
\newcommand{\cN}{\mathcal N}

\newcommand{\cP}{\mathcal P}

\newcommand{\cT}{\mathcal T}

\newcommand{\her}{{}_{\rm h}}

\newcommand{\id}{\mathds{1}} % identity matrix
\DeclareMathOperator{\ii}{i} % imaginary unit
\DeclareMathOperator*{\argmin}{argmin}
\DeclareMathOperator*{\argmax}{argmax}
\DeclareMathOperator{\conv}{conv}
\DeclareMathOperator{\diag}{diag}
\DeclareMathOperator{\tr}{tr}
\DeclareMathOperator{\pos}{pos}
\DeclareMathOperator{\cs}{cs}
\begin{document}
\selectlanguage{english}
\title{Operator systems and convex sets with many normal cones}
\author{Stephan Weis}
\begin{abstract}
The state space of an operator system of $n$-by-$n$ matrices has, in a 
sense, many normal cones. Merely this convex geometrical property implies
smoothness qualities and a clustering property of exposed faces. The 
latter holds since each exposed face is an intersection 
of maximal exposed faces. An isomorphism translates these results to the 
lattice of ground state projections of the operator system. We work on 
minimizing the assumptions under which a convex set has the mentioned 
properties. 
\end{abstract}
%primary
%52A20 Convex sets in $n$ dimensions
%52B05 Combinatorial properties (number of faces, shortest paths, etc.)
%47L07 Convex sets and cones of operators
%51D25 Lattices of subspaces
%47A12 Numerical range, numerical radius
%secundary
%81P16 Quantum state spaces, operational and probabilistic concepts
%14P10 Semialgebraic sets and related spaces
%
\maketitle
\thispagestyle{empty}
\pagestyle{myheadings}
\markleft{\hfill Operator systems and convex sets with many normal cones\hfill}
\markright{\hfill S.\ Weis \hfill}
%
%%%%%%%%%%%%%%%%%%%%%%%%%%%%%%%%%%%%%%%%%%%%%%%%%%%%%%%%%%%%%%%%%%%%%%%%%%%%
%%%%%%%%%%%%%%%%%%%%%%%%%%%%%%%%%%%%%%%%%%%%%%%%%%%%%%%%%%%%%%%%%%%%%%%%%%%%
%%%%%%%%%%%%%%%%%%%%%%%%%%%%%%%%%%%%%%%%%%%%%%%%%%%%%%%%%%%%%%%%%%%%%%%%%%%%
%%%%%%%%%%%%%%%%%%%%%%%%%%%%%%%%%%%%%%%%%%%%%%%%%%%%%%%%%%%%%%%%%%%%%%%%%%%%
%%%%%%%%%%%%%%%%%%%%%%%%%%%%%%%%%%%%%%%%%%%%%%%%%%%%%%%%%%%%%%%%%%%%%%%%%%%%
%
\section{Introduction}
\label{sec:intro}
\par
We are interested in the convex geometry of a reduced statistical model
of quantum mechanics. The set of density matrices $\cM_n$ consists of 
positive-semidefinite hermitian matrices of trace one and represents the 
{\em state space} \cite{AlfsenShultz2001} of normalized positive linear 
functionals on the C*-algebra $M_n$ of complex $n$-by-$n$ matrices, 
$n\in\bN$. This convex set is a statistical model of quantum mechanics 
\cite{Holevo2011,BengtssonZyczkowski2006}. 
Let $S\subset M_n$ be an {\em operator system} \cite{Paulsen2002}, that is
a complex vector space which contains the identity matrix $\id$ and which 
is self-adjoint, that is $s\in S$ implies $s^*\in S$. 
The state space of $S$, that is the set of normalized positive linear 
functionals on $S$, is represented by the projection $\cM(S)$ of $\cM_n$ 
onto $S$, which we call {\em state space} in the following. This 
reduction of $\cM_n$ has a very broad use in quantum mechanics, for 
example in state tomography \cite{Heinosaari-etal2013}, inference 
\cite{Jaynes1957}, or quantum chemistry \cite{Erdahl1972},
because $\cM(S)$ represents expected values of observables, probabilities 
of measurements (POVM's), and reduced density matrices (quantum marginals).
\par
Coordinate representations of $\cM(S)$ are known in operator theory as
the convex hull of the {\em joint numerical range} 
\cite{Au-YeungPoon1979,LiPoon2000,Gutkin-etal2004,KrupnikSpitkovsky2006,
ChienNakazato2010,Cheung-etal2011,GutkinZyczkowski2013,Szymanski-etal2016}
or {\em joint algebraic numerical range} \cite{Mueller2010}. They are 
{\em algebraic polars} of spectrahedra \cite{RamanaGoldman1995} which makes 
convex algebraic geometry 
\cite{Netzer-etal2010,RostalskiSturmfels2013,SinnSturmfels2015,Magron-etal2015}
useful to study $\cM(S)$. In analogy with statistics 
\cite{Barndorff-Nielsen1978} we called coordinate representations of 
$\cM(S)$ {\em convex support sets} \cite{Weis2011}. 
\par
Convex geometry of $\cM(S)$ is highlighted by signatures of quantum phase 
transitions that appear already for quantum systems with a finite 
dimensional state space $\cM_n$. These signatures are marked by abrupt 
changes of maximum-entropy states and ground state projections 
\cite{Arrachea-etal1992,Chen-etal2015}. These quantities have discontinuities 
which are indeed related to the convex geometry of $\cM(S)$ 
\cite{WeisKnauf2012,Weis2014,Chen-etal2015,Rodman-etal2016}.
The present article focusses on the notion of {\em exposed face} 
of a convex set $C\subset\bR^m$, which is either the empty set or the set of 
minimizers of a linear form on $C$. We denote the set of exposed faces of $C$ 
by $\cE_C$. On $\cE_C$ we consider the partial ordering by inclusion.
The general theory of the partial ordering of $\cE_C$ is well-understood 
\cite{Barker1978,LoewyTam1986}. 
\par
It is not surprising that $\cE_{\cM(S)}$ is order isomorphic to the lattice 
of ground state projections of hermitian matrices in $S$. Indeed, the ground 
state energy of a quantum mechanical operator, represented by a hermitian 
matrix $s$, is the least eigenvalue of $s$ and the minimum is achieved for 
density matrices which are concentrated on the corresponding eigenspace, the 
ground state space. The ground state projection is the orthogonal projection 
onto the ground state space. 
\par
Earlier work \cite{Chen-etal2012a} on the ordering of ground state spaces 
addresses the space $X$ of {\em $k$-local Hamiltonians}, which is frequently
used in many-party physics. The state space $C:=\cM(X+\ii X)$ represents the 
set of $k$-party reduced density matrices, whose convex geometry has a longer 
history \cite{Erdahl1972} and is still a topic 
\cite{Ocko-etal2011,Chen-etal2012b,Chen-etal2012c}. The partial ordering of 
$\cE_C$ was studied from the point of view of  minimal elements of 
$\cE_C\setminus\{\emptyset\}$ and their suprema \cite{Chen-etal2012a}. The 
present article continues our work \cite{Weis2012a} to study maximal elements 
of $\cE_C\setminus\{C\}$ and their infima for arbitrary convex sets $C$ (not 
necessarily closed or bounded). Properties of $\cM(S)$ are reflected in $C$
only by assumptions on normal cones.
%
%
%%%%%%%%%%%%%%%%%%%%%%%%%%%%%%%%%%%%%%%%%%%%%%%%%%%%%%%%%%%%%%%%%%%%%%%%%%%%
%%%%%%%%%%%%%%%%%%%%%%%%%%%%%%%%%%%%%%%%%%%%%%%%%%%%%%%%%%%%%%%%%%%%%%%%%%%%
%%%%%%%%%%%%%%%%%%%%%%%%%%%%%%%%%%%%%%%%%%%%%%%%%%%%%%%%%%%%%%%%%%%%%%%%%%%%
%%%%%%%%%%%%%%%%%%%%%%%%%%%%%%%%%%%%%%%%%%%%%%%%%%%%%%%%%%%%%%%%%%%%%%%%%%%%
%%%%%%%%%%%%%%%%%%%%%%%%%%%%%%%%%%%%%%%%%%%%%%%%%%%%%%%%%%%%%%%%%%%%%%%%%%%%
%
\section{Discussion of the Results}
\label{sec:results}
\par
Let $C\subset\bR^m$, $m\in\bN$, be a convex subset. The {\em normal cone} of 
$C$ at $x\in C$ is
\[
N_C(x):=\{u\in\bR^m\mid\forall y\in C:\,\langle u,y-x\rangle\leq 0\}
\]
where $\langle\cdot,\cdot\rangle$ is the standard inner product.
Elements of $N_C(x)$ are (outward pointing) {\em normal vectors} of $C$ 
at $x$. The {\em relative interior} of a convex subset $C\subset\bR^m$ 
is the interior of $C$ with respect to the topology of the affine 
hull of $C$. The {\em normal cone} $N_C(F)$ of $C$ at a non-empty 
convex subset $F$ of $C$ is well-defined as the normal 
cone $N_C(x)$ at any relative interior point $x$ of $F$ (see 
Section~2.2 of \cite{Schneider2014} or Section~4 of \cite{Weis2012a}). 
Let $\cN_C$ denote the set of normal cones at points of $C$ together
with $N_C(\emptyset):=\bR^m$. If $C$ is not a singleton then 
$N_C:\cE_C\to\cN_C$ is an antitone lattice isomorphism 
(\ref{eq:exp-norm-iso}).
\par
This article addresses exposed faces but the simple Example~\ref{exa:intro}(3) 
shows that state spaces do have a richer convex geometry.
A {\em face} of $C$ is a convex subset of $C$ containing every closed segment 
in $C$ whose relative interior it intersects. Exposed faces are faces, a
{\em non-exposed face} is a face which is not an exposed face. If a face $F$ 
is a singleton then its element is an {\em extreme point}, 
{\em exposed point}, or {\em non-exposed point} of $C$ respectively, if $F$ 
is a face, exposed face, or non-exposed face. We use the analogous definitions 
for faces which are rays.
\par
To exclude trivialities we call an exposed face $F\in\cE_C$ 
{\em proper exposed face} of $C$, if $F\not\in\{\emptyset,C\}$. We define a 
{\em proper convex subset} of $\bR^m$ to be a convex subset of $\bR^m$ with 
interior points\footnote{%
The simplifying assumption of interior points is fulfilled by any convex set 
after applying an affine embedding which removes codimensions. It guarantees 
that the normal cone at every point is a pointed closed convex cone.
By definition, a {\em convex cone} is a non-empty convex subset $K$ of 
$\bR^m$ such that $\alpha x\in K$ for all $\alpha\geq 0$ and $x\in K$. 
A convex cone $K$ is {\em pointed}, if $K\cap(-K)=\{0\}$.}
which has a proper exposed face. A {\em proper normal cone} of $C$ is any
element of $\cN_C$ other than $N_C(\emptyset)$ and $N_C(C)$. 
\begin{defn}\label{def:main}
Let $\cC_m^0$, $\cC_m$, $\cC_m'$, and $\cC_m''$ denote, respectively, the
class of proper convex subsets of $\bR^m$, $m\in\bN$, such that for every 
proper normal cone $N$ of $C$\\[.5\baselineskip]
\centerline{%
\def\arraystretch{2}
\begin{tabular}{l|l}
\hline
\parbox{1cm}{\centering $\cC_m^0$} & 
\parbox{13cm}{%
$N$ has an exposed ray which is in $\cN_C$,}\\
\hline
\parbox{1cm}{\centering $\cC_m$} & 
\parbox[t]{13cm}{%
$N$ has $\dim(N)$ linearly independent exposed rays which are 
in $\cN_C$,}\\
\hline
\parbox{1cm}{\centering $\cC_m'$} & 
\parbox[t]{13cm}{%
every extreme ray of $N$ is in $\cN_C$,}\\
\hline
\parbox{1cm}{\centering $\cC_m''$} & 
\parbox[t]{13cm}{%
every non-empty face of $N$ is in $\cN_C$.}\\
\hline
\end{tabular}}
\end{defn}
\vspace{.2mm}
\par
Let us critically discuss this definition.
We call {\em convex body} a compact convex set.
\begin{rem}\label{rem:def}~
\begin{enumerate}
\item
Replacing {\em exposed ray} with {\em extreme ray} does not change 
the definition of $\cC_m^0$ or $\cC_m$ because normal cones of $C$ 
included in $N$ are exposed faces of $N$ 
(Lemma~\ref{lem:normal-cone-inclucions-exposed}). 
\item
Replacing {\em extreme ray} with {\em exposed ray} weakens the definition 
of $\cC_m'$ (Example~\ref{exa:AltCm1}). A convex body lies in $\cC_m'$ if 
and only if its polar convex body 
has no non-exposed points (Theorem~\ref{thm:def-main}). 
Corollary~\ref{cor:coorbitope}, a sufficient condition for inclusion to
$\cC_m'$, shows that coorbitopes\footnote{% 
A {\em coorbitope} \cite{Sanyal-etal2011} is polar to an orbitope, where 
an {\em orbitope} is defined as the convex hull of the orbit of a compact 
algebraic group acting linearly on a vector space. Proposition~2.2 
of \cite{Sanyal-etal2011} shows that orbitopes have no non-exposed points 
because their extreme points lie on a sphere.}
lie in $\cC_m'$.
\item
Replacing {\em face} with {\em exposed face} does not change the definition 
of $\cC_m''$ (Lemma~\ref{lem:AltCm2}). A convex body lies in $\cC_m''$ if 
and only if its polar convex body has no non-exposed faces
(Theorem~\ref{thm:def-main}). In particular, state spaces of operator systems
lie in $\cC_m''$ (Corollary~\ref{cor:state-spaces-C2}).
\item
Lemma~\ref{lem:max-number-rays} proves $\cC_m' \subset \cC_m$ and this 
implies $\cC_m'' \subset \cC_m' \subset \cC_m \subset \cC_m^0$. See 
(\ref{eq:nesting}) for a refinement of this nesting of classes and for 
examples of strictness.
\end{enumerate}
\end{rem}
\par
Section~\ref{sec:smoothness} studies the class $\cC_m^0$ with a focus
on smoothness. We call {\em coatom}\footnote{%
Note that a coatom may not be a {\em facet}, that is a face of codimension 
one \cite{Schneider2014,Ziegler1995}. For a polyhedral convex set the 
notions of coatom and facet are equivalent \cite{Gruenbaum2003,Ziegler1995}.}
of $\cE_C$ an inclusion maximal element of $\cE_C\setminus\{C\}$. 
Theorem~\ref{thm:one-ray} proves that a proper convex subset $C$ of $\bR^m$ 
lies in $\cC_m^0$ if and only if every coatom of $\cE_C$ is a smooth exposed
face (the converse is trivial). Thereby an exposed face is smooth if it has a 
unique unit normal vector. Theorem~\ref{thm:char-C0} proves that $C\in\cC_m^0$ 
is equivalent to the boundary of $C$ being covered by smooth coatoms of 
$\cE_C$. See Example~\ref{exa:intro}(1) for a convex set without this property.
\par
Section~\ref{sec:proof-main} improves the theorem in Section~1.2.4 of 
\cite{Weis2012a}, which shows for $C\in\cC_m''$ that every proper exposed 
face of $C$ is an intersection of coatoms of $\cE_C$. This property is 
well-known for polytopes \cite{Ziegler1995}, which are included in $\cC_m''$ 
because they are convex support sets, see 
Remark~\ref{rem:convex-support} and 
Corollary~\ref{cor:state-spaces-C2}.
Theorem~\ref{thm:max-number-coatoms} weakens not only the assumptions from 
$C\in\cC_m''$ to $C\in\cC_m$ but adds a dimension dependent multiplicity:
\begin{cor}[Intersections]\label{cor:main}
Let $C\in\cC_m$ and let $F\in\cE_C$ be a proper exposed face. Then there 
exist $\dim(N_C(F))$ mutually distinct coatoms of $\cE_C$ whose intersection
is $F$ and whose normal cones are linearly independent exposed rays of 
$N_C(F)$.
\end{cor}
\noindent
Of course, if $\dim(N_C(F))>2$ then $F$ can be the intersection of any 
number (at least two) of coatoms. An example is an octahedron where 
every vertex is the intersection of two, three, or four facets. Another 
example is a cone based on a disk whose apex is the intersection of two 
surface lines while lying on a continuum of them.
\par
Corollary~\ref{cor:main} contains a method to construct clusters of 
exposed faces of $C\in\cC_m$. More precisely, we define a {\em cluster} 
as an equivalence class of coatoms of $\cE_C$ where two coatoms 
$F_1,F_2$ are equivalent if there is a sequence of coatoms $G_0,\ldots,G_k$, 
$k\in\bN$, such that $G_0=F_1$, $G_k=F_2$ and 
$G_{i-1}\cap G_i\neq\emptyset$ for $i=1,\ldots,k$. 
\begin{cor}[Construction of clusters]\label{cor:clustering}
Let $C\in\cC_m$ and let $F\in\cE_C$ be a proper exposed face. If $F$ 
strictly contains a proper exposed face $G\in\cE_C$ then there exists a 
proper exposed face $F'\neq F$ of $C$ such that $G\subset F'$.
One can choose $F'$ to be a coatom of $\cE_C$.
\end{cor}
\par
The described corollaries may be checked with 3D sets of the last 
paragraph of Example~\ref{exa:AltCm1} or Example~\ref{exa:JNR3x3}. 
2D examples suffice to distinguish the classes $\cC_m\subset\cC_m^0$. 
\begin{exa}\label{exa:intro}
Let $\conv(X)$ denote the convex  hull of a subset $X\subset\bR^m$. We study
\begin{enumerate}
\item
the lens
$C_{\rm ()}:=\{(x,y)\in\bR^2\mid (x\pm\tfrac{3}{2})^2+y^2\leq(\tfrac{5}{2})^2\}$,
\item
the truncated disk
$C_{\rm (|}:=\{(x,y)\in\bR^2\mid x^2+y^2\leq 1, y\leq\tfrac{1}{2}\}$,
\item
and the drop
$C_{\rm 0>}:=\conv(\{(x,y)\in\bR^2\mid x^2+y^2=1\}\cup\{a\})$,
\qquad $a:=(0,2)$.
\end{enumerate}
\par
We have $C_{\rm ()}\not\in\cC_2^0$ because the coatom $\{a\}$ of 
$\cE_{C_{\rm ()}}$ has a 2D normal cone. The exposed point 
$t_\pm:=\tfrac{1}{2}(\pm\sqrt{3},1)$ of $C_{\rm (|}$ has a 2D normal 
cone, too, but this normal cone has an exposed ray which is a normal 
cone of $C$, so $C_{\rm (|}\in\cC_2^0$ holds. Since this cone has 
only one exposed ray which is a normal cone of $C$ we have 
$C_{\rm (|}\in\cC_2^0\setminus\cC_2$. Another argument for 
$C_{\rm (|}\not\in\cC_2$ is that the construction of clusters of 
Corollary~\ref{cor:clustering} fails because $t_\pm$ lies only on one 
coatom of $\cE_C$, which is the segment $[t_-,t_+]$. 
\par
The drop $C_{\rm 0>}$ belongs to $\cC_2$ because $a$ is the only point 
with a 2D normal cone and because the exposed rays of this normal cone 
are normal cones. Alternatively, $C_{\rm 0>}\in\cC_2$ by 
Theorem~\ref{thm:def-main}(4) since $C_{\rm 0>}$ is the polar of 
$C_{\rm (|}$, which has no non-exposed faces. The inclusion 
$C_{\rm 0>}\in\cC_2$ follows also from Remark~\ref{rem:convex-support}(1)
and Corollary~\ref{cor:state-spaces-C2} since $C_{\rm 0>}$ is the convex 
support of 
$F_1:=\left(\begin{smallmatrix}0&1&0\\1&0&0\\0&0&0\end{smallmatrix}\right)$
and
$F_2:=\left(\begin{smallmatrix}0&-\ii&0\\\ii&0&0\\0&0&2\end{smallmatrix}\right)$. 
Notice that the extreme point $t_\pm$ of $C_{\rm 0>}$ is a non-exposed point. 
Since $C_{\rm 0>}\in\cC_2$ holds, this follows also by contradiction from
Corollary~\ref{cor:clustering}: If $t_\pm$ was an exposed point of $C_{\rm 0>}$ 
then it had to lie on two coatoms of $\cE_{C_{\rm 0>}}$.
\end{exa}
\par
Dimension three is needed to differentiate between the classes 
$\cC_m'' \subset \cC_m' \subset \cC_m$ (because $\cC_2''=\cC_2$ holds)
and to finish Remark~\ref{rem:def}. This discussion will be done
with convex bodies for which Theorem~\ref{thm:def-main} translates
Definition~\ref{def:main} to polar convex bodies. Further, 
Section~\ref{sec:convex-bodies} recalls that the class $\cC_m''$ is 
closed under projection to subspaces, which is wrong for 
$\cC_m^0,\cC_m$, and $\cC_m'$.
\par
Theorem~\ref{thm:state-space} proves that the state space $\cM(S)$ of 
an operator system $S\subset M_n$ is the projection of the state space 
of the algebra $M_n$ onto $S$. Further topics of 
Section~\ref{sec:state-spaces} are a proof of 
$\cM(S)\in\cC_{\dim_\bC(S)-1}''$, coordinate representations of 
$\cM(S)$ in terms of convex support sets and joint numerical ranges, 
the isomorphism between exposed faces of $\cM(S)$ and ground state 
projections, and a discussion of state spaces of $3$-by-$3$ matrices. 
%
%%%%%%%%%%%%%%%%%%%%%%%%%%%%%%%%%%%%%%%%%%%%%%%%%%%%%%%%%%%%%%%%%%%%%%%%%%%%
%%%%%%%%%%%%%%%%%%%%%%%%%%%%%%%%%%%%%%%%%%%%%%%%%%%%%%%%%%%%%%%%%%%%%%%%%%%%
%%%%%%%%%%%%%%%%%%%%%%%%%%%%%%%%%%%%%%%%%%%%%%%%%%%%%%%%%%%%%%%%%%%%%%%%%%%%
%%%%%%%%%%%%%%%%%%%%%%%%%%%%%%%%%%%%%%%%%%%%%%%%%%%%%%%%%%%%%%%%%%%%%%%%%%%%
%%%%%%%%%%%%%%%%%%%%%%%%%%%%%%%%%%%%%%%%%%%%%%%%%%%%%%%%%%%%%%%%%%%%%%%%%%%%
%
\section{Smoothness}
\label{sec:smoothness}
\par
We characterize the class $\cC_m^0$ in terms of smoothness properties.
Further, we study some very special convex sets which are smooth,
strictly convex, or both (ovals).
\par
A {\em lattice} is a partially ordered set where the infimum and 
supremum of each two elements exists. An {\em atom} in a lattice 
$(\cL,\leq,\wedge,\vee,0)$ with smallest element $0$ is an element 
$x\in\cL$, $x\neq 0$, such that $y\leq x$, $y\neq x$ implies $y=0$ 
for all $y\in\cL$. Similarly, a {\em coatom} in a lattice 
$(\cL,\leq,\wedge,\vee,1)$ with greatest element $1$ is an element 
$x\in\cL$, $x\neq 1$, such that $y\geq x$, $y\neq x$ implies $y=1$ 
for all $y\in\cL$. The lattice $(\cL,\leq,\wedge,\vee,0)$ is 
{\em atomistic} if every element of $\cL$ is the supremum of the atoms 
which it contains \cite{LoewyTam1986} (such a lattice is called 
{\em atomic} in \cite{Birkhoff1967,Ziegler1995}). The lattice
$(\cL,\leq,\wedge,\vee,1)$ is {\em coatomistic} if every element is the 
infimum of the coatoms in which it is contained. A lattice is 
{\em complete} if an arbitrary subset has an infimum and a supremum.
The smallest and greatest elements of a lattice, when they exist, 
are called {\em improper} elements, all other elements are {\em proper} 
elements. If $x,y\in\cL$ and $x\leq y$ then we define the {\em interval}
$[x,y]_\cL:=\{z\in\cL\mid x\leq z\leq y\}$.
\par
Without reminder we will use the fact that the smallest exposed face 
containing a proper face is a proper exposed face (Lemma 4.6 of 
\cite{Weis2012a}). In particular, for a convex subset of $\bR^m$ to have 
a proper exposed face is the same as to have a proper face.
\par
Let $C\subset\bR^m$, $m\in\bN$, be a convex subset. Both $\cE_C$ 
and $\cN_C$ are, partially ordered by inclusion, complete lattices where 
the infimum is the intersection \cite{Weis2012a}. The improper elements 
of $\cE_C$ are $\emptyset$ and $C$ and the improper elements of $\cN_C$ 
are $N_C(\emptyset)=\bE^m$ and $N_C(C)$, the latter being the vector 
space which is the orthogonal complement of the affine hull of $C$. By 
Proposition~4.7 of \cite{Weis2012a}, if $C$ is not a singleton then
\begin{equation}\label{eq:exp-norm-iso}
N_C:\cE_C\to\cN_C
\end{equation}
is an antitone lattice isomorphism. Let $C$ be a proper convex subset of 
$\bR^m$. Then the coatoms of $\cE_C$ are proper exposed faces 
and the atoms of $\cN_C$ proper normal cones. Moreover, proper normal 
cones are pointed closed convex cones and the smallest element of $\cN_C$ 
is $\{0\}$.
\par
We start with an observation about faces of normal cones. 
The {\em positive hull} of a non-empty subset $X\subset\bR^m$ is 
$\pos(X):=\{\alpha x\mid\alpha\geq 0, x\in X\}$.
\begin{lem}\label{lem:normal-cone-inclucions-exposed}
Let $C\subset\bR^m$, $m\in\bN$, be a convex set. If $K,L\in\cN_C$,
$K\subset L$, and $L$ is a proper normal cone of $C$ then $K$ is 
an exposed face of $L$.
\end{lem}
{\em Proof:} 
Using the antitone lattice isomorphism (\ref{eq:exp-norm-iso}) and 
applying a translation to $C$ we can assume that $0\in C$, and that 
$K=N_C(x)$ and $L=N_C(0)$ hold for some $x\in C$. The equalities
\[
N_C(0)=N_{\pos(C)}(0),
\qquad
N_C(x)=N_{\pos(C)}(x)
\]
are easy to prove, the second equality holds because all points $y$ 
of the segment $[x,0[$ have the same normal cone $N_C(y)=N_C(x)$.
Secondly, it is easy to see that $N_{\pos(C)}(x)$ is the intersection 
of $N_{\pos(C)}(0)$ with the orthogonal complement of the face of 
$\pos(C)$ containing $x$ in its relative interior. This intersection
is an exposed face of $N_{\pos(C)}(0)$, known as {\em dual face}
\cite{Tam1985}.
\hspace*{\fill}$\square$\\
\par
Let $C\subset\bR^m$ be a convex subset with interior points. A point
$s\in C$ is a {\em smooth point} \cite{Schneider2014} of $C$ if
$N_C(s)$ is one-dimensional. In that case 
$N_C(s)$ is a ray and an atom of $\cN_C$. Similarly we call 
$F\in\cE_C$ a {\em smooth exposed face} of $C$ if $N_C(F)$ is a ray
and we call $C$ a {\em smooth convex set} if all proper normal cones
of $C$ are rays. 
\begin{thm}\label{thm:one-ray}
Let $C$ be a proper convex subset of $\bR^m$, $m\in\bN$. Then $C\in\cC_m^0$ 
holds if and only if every atom of $\cN_C$ is a ray. In that case the 
isomorphism (\ref{eq:exp-norm-iso}) restricts to a bijection from the 
coatoms of $\cE_C$ to the normal cones of $C$ which are rays.
\end{thm}
{\em Proof:} 
Let $C\in\cC_m^0$ and let $N$ be an atom of $\cN_C$. Then $N$ is a proper
normal cone and by definition of $\cC_m^0$ the pointed cone $N$ has an
exposed ray which is a normal cone of $C$. Since $N$ is minimal in
$\cN_C\setminus\{\{0\}\}$, it must be equal to that ray.
Conversely, we recall the well-known property that the face lattice of a
finite-dimensional convex set has finite length, so $\cE_C$ has finite
length and by the isomorphism (\ref{eq:exp-norm-iso}) $\cN_C$ has finite 
length. It follows that every proper normal cone $N$ of $C$ contains an 
atom $K$ of $\cN_C$. By assumption, $K$ is a ray and by 
Lemma~\ref{lem:normal-cone-inclucions-exposed} the ray $K$ is an exposed
ray of $N$. This proves the first statement. The second statement 
follows from the first one because the normal cones which are rays are 
atoms of $\cN_C$.
\hspace*{\fill}$\square$\\
\par
A {\em sublattice} \cite{Birkhoff1967} of a lattice 
$(\cL,\leq,\wedge,\vee)$ is a subset $\cL'\subset\cL$ such that 
$a,b\in\cL'$ implies that $a\wedge b$ and $a\vee b$ lie in $\cL'$. 
Clearly, the intersection of any family of sublattices is another 
sublattice. Therefore the smallest sublattice of $\cL$ containing 
a given subset $X\subset\cL$ exists, we call it the sublattice 
{\em generated} by $X$. 
\par
Let $C\subset\bR^m$ be a convex subset with interior points.
Every ray in $\cN_C$ is an atom of $\cN_C$ so every smooth exposed 
face of $C$ is a coatom of $\cE_C$ by (\ref{eq:exp-norm-iso}). 
Hence, every smooth exposed point $s$ of $C$ is simultaneously a 
coatom and an atom of $\cE_C$ (we identify a point with its 
singleton). Therefore $\{\emptyset,\{s\},C\}$ is a sublattice of 
$\cE_C$ and the remainder $\cE_C\setminus\{s\}$ is a sublattice, 
too. Similarly, if we define
\[
S_C:=\{s\in C\mid\mbox{$s$ is a smooth exposed point of $C$}\}
\]
then $S_C\cup\{\emptyset,C\}$ and $\cE_C\setminus S_C$ are sublattices 
of $\cE_C$. Let $\cL_C$ be the sublattice of $\cE_C$ generated by 
$\{\emptyset,C\}$ and by the smooth coatoms of $\cE_C$ which are no 
singletons. We call $\cL_C$ {\em lattice of large smooth coatoms} of 
$\cE_C$. By what we have discussed in this paragraph, 
$\cL_C\subset\cE_C\setminus S_C$ holds.
\begin{thm}\label{thm:char-C0}
Let $C$ be a proper convex subset of $\bR^m$, $m\in\bN$. Then $S_C$ and 
$\bigcup_{F\in\cL_C\setminus C}F$ are disjoint. We have $C\in\cC_m^0$ 
if and only if the boundary of $C$ is the union of the smooth coatoms 
of $\cE_C$.
\end{thm}
{\em Proof:} 
About the first statement, we have noticed in the preceding paragraph
that $\cL_C\subset\cE_C\setminus S_C$ holds, so for every $s\in S_C$ 
it follow $\{s\}\not\in\cL_C$. Since $\{s\}$ is a coatom of $\cE_C$ 
the point $s$ cannot lie on any coatom of $\cE_C$ other than $\{s\}$. 
This proves the disjointness.
\par
We derive the second statement from the property that the union of coatoms 
of $\cE_C$ is the boundary of $C$, which holds because the union of 
proper exposed faces of $C$ is the boundary of $C$ (see for example Theorem 
18.2 of~\cite{Rockafellar1970}). If $C\in\cC_m^0$ then 
Theorem~\ref{thm:one-ray} shows that every coatom of $\cE_C$ is smooth. The 
converse is true because every coatom of $\cE_C$ is needed to fill up the 
boundary of $C$. Indeed, if a relative interior point $x$ of a coatom $F$ 
of $C$ is covered by a coatom $G$ of $\cE$ then $G\subset F$ follows from 
Theorem~18.1 of \cite{Rockafellar1970}. This proves $F=G$.
\hspace*{\fill}$\square$\\
\par
Example~\ref{exa:JNR3x3} discusses the partition of the boundary of
$C$ into $S_C$ and its complement for convex support sets.
Notice that the second assertion of Theorem~\ref{thm:char-C0} is 
weaker than $\cE_C=S_C\cup\cL_C$, which holds for 
$C\in\cC_m$ by Theorem~\ref{thm:max-number-coatoms}. The assumptions
of these two theorems cannot be weakened arbitrarily: The boundary of 
the lens $C_{\rm ()}\not\in\cC_2^0$, see Example~\ref{exa:intro}, 
is not covered by smooth coatoms of $\cE_C$. We have a proper 
inclusion $S_C\cup\cL_C\subsetneq\cE_C$ for the truncated disk 
$C=C_{\rm (|}\in\cC_2^0\setminus\cC_2$.
\par
We now collect properties of strictly convex bodies, which we meet
again in Section~\ref{sec:convex-bodies} as the polars of smooth 
convex bodies. Let $C\subset\bR^m$ be a convex subset with interior 
points. We call $C$ {\em strictly convex} if the relative interior
of every closed segment in $C$ lies in the interior of $C$, 
{\em pseudo-oval} if C is smooth and strictly convex, and {\em oval} 
if $C$ is a compact pseudo-oval.
\begin{lem}\label{lem:char-strict}
Let $C$ be a proper convex subset of $\bR^m$, $m\in\bN$. Then the 
following statements are equivalent.
\begin{enumerate}
\item
The set $C$ is strictly convex,
\item
every coatom of $\cE_C$ is a singleton,
\item
every boundary point of $C$ is an exposed point of $C$.
\end{enumerate}
\end{lem}
{\em Proof:} 
We show (1)$\implies$(2)$\implies$(3)$\implies$(1).
The first implication is true because coatoms of $\cE_C$ are 
non-empty convex subsets of the boundary of $C$, see the second 
paragraph of the proof of Theorem~\ref{thm:char-C0}. The second 
implication is true because every boundary point of $C$ lies in 
a proper exposed face of $C$. The third implication follows from 
the definition of a face.
\hspace*{\fill}$\square$\\
\par
Lemma~\ref{lem:char-strict} and the following statement improve 
Theorem~3.1 of \cite{Szymanski-etal2016} about ovals. Notice
that the lens $C_{\rm ()}\not\in\cC_2^0$ from Example~\ref{exa:intro} 
is strictly convex but no pseudo-oval.
\begin{cor}\label{cor:char-ovals}
Let $C\in\cC_m^0$, $m\in\bN$. If $C$ is strictly convex then $C$ 
is a pseudo-oval.
\end{cor}
{\em Proof:}
Theorem~\ref{thm:char-C0} shows that the boundary of 
$C$ is covered by smooth coatoms of $\cE_C$, which are singletons by 
Lemma~\ref{lem:char-strict} because $C$ is strictly convex. Therefore 
every proper normal cone of $C$ is a ray which proves the claim.
\hspace*{\fill}$\square$\\
\par
The next statement characterizes ovals under a stronger assumption 
than required in Corollary~\ref{cor:char-ovals}. The proof uses 
{\em Minkowski's theorem}, which states that every non-empty convex 
body is the convex hull of its extreme points (for a proof see for 
example Corollary 1.4.5 of \cite{Schneider2014}). 
\begin{lem}\label{lem:smooth-oval}
Let $C$ be a proper smooth convex body in $\bR^m$, $m\in\bN$.
Then $C$ is an oval if and only if $C$ has no non-exposed face.
\end{lem}
{\em Proof:} 
An oval has no non-exposed face because this would imply the existence 
of a proper exposed face of dimension one or larger. Conversely, by
Lemma~\ref{lem:char-strict} it suffices to show that any coatom 
$F$ of $\cE_C$ is a singleton. Minkowski's theorem shows that $F$ 
contains an extreme point $x$ which, by hypothesis, is an exposed 
point of $C$. By contradiction, if $\{x\}\neq F$ then the isomorphism 
(\ref{eq:exp-norm-iso}) shows that the normal cone $N_C(x)$ has at least 
dimension two, which contracts the smoothness of $C$.
\hspace*{\fill}$\square$\\
\par
Analogues of Lemma~\ref{lem:smooth-oval} about unbounded and non-closed 
sets are wrong because of possibly  missing extreme points. Consider a 
closed half-space or an open square with an open segment attached to one 
of its sides.
%
%%%%%%%%%%%%%%%%%%%%%%%%%%%%%%%%%%%%%%%%%%%%%%%%%%%%%%%%%%%%%%%%%%%%%%%%%%%%
%%%%%%%%%%%%%%%%%%%%%%%%%%%%%%%%%%%%%%%%%%%%%%%%%%%%%%%%%%%%%%%%%%%%%%%%%%%%
%%%%%%%%%%%%%%%%%%%%%%%%%%%%%%%%%%%%%%%%%%%%%%%%%%%%%%%%%%%%%%%%%%%%%%%%%%%%
%%%%%%%%%%%%%%%%%%%%%%%%%%%%%%%%%%%%%%%%%%%%%%%%%%%%%%%%%%%%%%%%%%%%%%%%%%%%
%%%%%%%%%%%%%%%%%%%%%%%%%%%%%%%%%%%%%%%%%%%%%%%%%%%%%%%%%%%%%%%%%%%%%%%%%%%%
%
\section{Intersections of exposed faces}
\label{sec:proof-main}
\par
We show that every proper exposed face of a set of class $\cC_m$ admits a 
representation as an intersection of coatoms of exposed faces, taking
into account the dimension of normal cones. At the end of the section we 
continue the discussion of the classes $\cC_m''\subset\cC_m'\subset\cC_m$ 
from Remark~\ref{rem:def}.
\par
Corollary~\ref{cor:main} is a consequence of the following statement. 
\begin{thm}\label{thm:max-number-coatoms}
Let $C$ be a proper convex subset of $\bR^m$, $m\in\bN$, and let $F\in\cE_C$ 
be a proper exposed face, $N=N_C(F)$, and $d:=\dim(N)$. Then the following 
statements are equivalent.
\begin{enumerate}
\item
$N$ has $d$ linearly independent exposed rays which are normal 
cones of $C$. 
\item
there exist $d$ mutually distinct coatoms of $\cE_C$ whose intersection
is $F$ and whose normal cones are linearly independent exposed rays of $N$.
\end{enumerate}
\end{thm}
{\em Proof:} 
The statement (2) is clearly stronger than (1), we prove that it follows 
from (1). Let $N$ have $d$ linearly independent exposed rays which lie in 
$\cN_C$. Their supremum $K$ in $\cN_C$ is included in $N$, hence $K$ is 
an exposed face of $N$ by Lemma~\ref{lem:normal-cone-inclucions-exposed}. 
This shows $K=N$ because 
$\dim(K)=d$ and because a proper face of any convex set has 
codimension at least one (Corollary 8.1.3 of \cite{Rockafellar1970}). 
The rays are atoms of $\cN_C$, hence the isomorphism
(\ref{eq:exp-norm-iso}) shows that the corresponding exposed faces are 
coatoms of $\cE_C$ and that their intersection is $F$.
\hspace*{\fill}$\square$\\
\par
We discuss further corollaries of Theorem~\ref{thm:max-number-coatoms}.
Let $C\subset\bR^m$, $m\in\bN$, be any convex subset. We call a face 
$F\neq\emptyset$ of $C$ a {\em corner}, if $\dim(N_C(F))=m$. Every corner 
is an exposed face of $C$ and a singleton, see for example Lemma~4.4 of 
\cite{Weis2012b}; we call its element a {\em corner point}. 
\begin{cor}\label{cor:max-number-coatoms}
Let $C$ be a proper convex subset of $\bR^m$, $m\in\bN$, and let $F\in\cE_C$ 
be a proper exposed face for which the equivalent statements of 
Theorem~\ref{thm:max-number-coatoms} hold. Then the following is true.
\begin{enumerate}
\item
The normal cone $N_C(F)$ is a ray if and only if $F$ is a coatom of $\cE_C$,
\item
$\dim(N_C(F))=2$ holds if and only if $F$ lies on a unique pair of coatoms 
of $\cE_C$; in that case $F$ is the intersection of the pair,
\item
$F$ is a corner of $C$ if and only if $F$ is the intersection of $m$ 
mutually distinct coatoms of $\cE_C$ whose normal cones are exposed 
rays of $N_C(F)$ which span $\bR^m$.
\item
The lattice $[F,C]_{\cE_C}$ is coatomistic. The coatoms of 
$[F,C]_{\cE_C}$ are the coatoms of $\cE_C$ which are included in 
$[F,C]_{\cE_C}$.
\end{enumerate}
\end{cor}
{\em Proof:} 
This corollary follows mainly from Theorem~\ref{thm:max-number-coatoms}(2).
A case selection of the dimension if $N_C(F)$ suffices to complete the 
proofs of the Items (1)--(3). Item (4) is completed by using basic 
properties of intervals and coatoms.
\hspace*{\fill}$\square$\\
\par
Corollary~\ref{cor:max-number-coatoms}(4) has an interpretation in terms
of normal cones.
\begin{rem}
If $F$ is a proper exposed face of a proper convex subset $C\subset\bR^m$ 
then the isomorphism  (\ref{eq:exp-norm-iso}) restricts to an antitone 
lattice isomorphism 
\[
[F,C]_{\cE_C}\rightarrow [\{0\},N_C(F)]_{\cN_C}
\]
to the interval of normal cones $[\{0\},N_C(F)]_{\cN_C}$, which is therefore 
atomistic. This isomorphism is simplest possible for $C\in\cC_m''$. Then all 
non-empty faces of $N_C(F)$ lie in $\cN_C$ by 
Lemma~\ref{lem:normal-cone-inclucions-exposed} and they all belong to the 
interval $[\{0\},N_C(F)]_{\cN_C}$.
\end{rem}
\par
We have to clarify a subtlety when stating that $\cE_C$ is coatomistic.
Namely, the intersection of proper exposed faces may be non-empty for 
some non-closed sets and is indeed non-empty for all closed convex 
cones%
\footnote{If the intersection of coatoms of $\cE_C$ is non-empty then
$\cE_C\setminus\{\emptyset\}$ would be a more natural definition of a 
lattice of exposed faces of $C$.}. This is not so for convex bodies.
\begin{lem}\label{lem:coatom-intersection}
Let $C$ be a proper convex body of $\bR^m$, $m\in\bN$.
Then the intersection of coatoms of $\cE_C$ is empty.
\end{lem}
{\em Proof:} 
If the intersection of coatoms of $\cE_C$ is non-empty then by
(\ref{eq:exp-norm-iso}) there is a proper normal cone $N\in\cN_C$ 
which contains all proper normal cones. Since $N$ is pointed,
some vectors of $\bE^m$ are no normal vectors at points of $C$. 
This can only happen when $C$ is unbounded or not closed.
\hspace*{\fill}$\square$\\
\par
Corollary~\ref{cor:max-number-coatoms}(4) has two global formulations.
Lemma~\ref{lem:coatom-intersection} shows that
Corollary~\ref{cor:add}(1) applies to convex bodies.
\begin{cor}\label{cor:add}
Let $C\in\cC_m$, $m\in\bN$.
\begin{enumerate}
\item
If the intersection of coatoms of $\cE_C$ is empty then $\cE_C$ is 
coatomistic.
\item
If the intersection of coatoms of $\cE_C$ is $A\neq\emptyset$ then
$\cE_C\setminus\{\emptyset\}=[A,C]_{\cE_C}$ is a coatomistic lattice 
where the infimum is the intersection.
\end{enumerate}
\end{cor}
{\em Proof:} 
This follows immediately from Theorem~\ref{thm:max-number-coatoms}(2).
\hspace*{\fill}$\square$\\
\par
We show $\cC_m'\subset\cC_m$ using Straszewicz's theorem, 
which affirms that every extreme point of a convex body is a limit of 
exposed points (for a proof see for example Theorem 1.4.7 of 
\cite{Schneider2014}).
\begin{lem}\label{lem:max-number-rays}
Let $C$ be a proper convex subset of $\bR^m$, $m\in\bN$. Let $N\in\cN_C$ 
be a proper normal cone all exposed rays of which lie in $\cN_C$. Then 
$N$ has $\dim(N)$ linearly independent exposed rays which lie in $\cN_C$.
\end{lem}
{\em Proof:} 
Since $N$ is a closed pointed cone, it admits a compact hyperplane
intersection through its interior. Straszewicz's theorem shows that 
this intersection is the closed convex 
hull of its exposed points, so $N$ is the closed convex hull of its 
exposed rays. Hence $N$ has $\dim(N)$ exposed rays which are linearly 
independent. By assumption these rays are normal cones of $C$.
\hspace*{\fill}$\square$\\
\par
Actually, Lemma~\ref{lem:max-number-rays} proves a bit more than 
$\cC_m'\subset\cC_m$, we will return to it in (\ref{eq:nesting}).
We prove the claim of Remark~\ref{rem:def}(3) that the class 
$\cC_m''$ does not increase when {\em face} is replaced with 
{\em exposed face} in the Definition~\ref{def:main} of $\cC_m''$. 
\begin{defn}
Let $\cC_m''^*$ denote the class of proper convex subsets of $\bR^m$, 
$m\in\bN$, such that every non-empty exposed face of every proper 
normal cone of $C$ is a normal cone of $C$.
\end{defn}
\par
Let $C$ be a convex subset of $\bR^m$, $m\in\bN$. A subset $F\subset C$ 
is a {\em poonem} \cite{Gruenbaum2003} of $C$ if there are sets 
$F_0\subset\cdots\subset F_k$, $k\in\bN$, such that $F_0=F$ 
and $F_k=C$ and such that $F_{i-1}$ is an exposed face of $F_i$ for
$i=1,\ldots,k$. It is not hard to prove that $F\subset C$ is a poonem
of $C$ if and only if $F$ is a face of $C$, see Section 1.2.1 of
\cite{Weis2012a}.
\begin{lem}\label{lem:AltCm2}
For all $m\in\bN$ we have $\cC_m''=\cC_m''^*$.
\end{lem}
{\em Proof:} 
Clearly $\cC_m''\subset\cC_m''^*$ holds. Conversely let 
$C\in\cC_m''^*$. Because of the isomorphism (\ref{eq:exp-norm-iso})
between exposed faces and normal cones of $C$, it suffices to prove 
for every non-empty face $T$ of every proper normal cone $N$ of $C$
that $T$ is a normal cone of $C$. Since $T$ is a poonem of $N$ 
there are sets $F_0\subset\cdots\subset F_k$, $k\in\bN$, such that 
$F_0=T$ and $F_k=N$ and such that $F_{i-1}$ is an exposed face of 
$F_i$ for $i=1,\ldots,k$. By the assumption of $C\in\cC_m''^*$ the
exposed face $F_{k-1}$ of $N$ is a normal cone of $C$. By
induction $T$ is a normal cone of $C$.
\hspace*{\fill}$\square$\\
%
%%%%%%%%%%%%%%%%%%%%%%%%%%%%%%%%%%%%%%%%%%%%%%%%%%%%%%%%%%%%%%%%%%%%%%%%%%%%
%%%%%%%%%%%%%%%%%%%%%%%%%%%%%%%%%%%%%%%%%%%%%%%%%%%%%%%%%%%%%%%%%%%%%%%%%%%%
%%%%%%%%%%%%%%%%%%%%%%%%%%%%%%%%%%%%%%%%%%%%%%%%%%%%%%%%%%%%%%%%%%%%%%%%%%%%
%%%%%%%%%%%%%%%%%%%%%%%%%%%%%%%%%%%%%%%%%%%%%%%%%%%%%%%%%%%%%%%%%%%%%%%%%%%%
%%%%%%%%%%%%%%%%%%%%%%%%%%%%%%%%%%%%%%%%%%%%%%%%%%%%%%%%%%%%%%%%%%%%%%%%%%%%
%
\section{Convex bodies}
\label{sec:convex-bodies}
\par
Theorem~\ref{thm:def-main} is an equivalent statement of 
Definition~\ref{def:main} for convex bodies in terms of their polars. 
It allows to present the Examples~\ref{exa:c1nonc2}, \ref{exa:AltCm1}, 
\ref{exa:cnonc1} which finish Remark~\ref{rem:def}. 
Remark~\ref{rem:projection} recalls that $\cC_m''$ is closed under 
projection to subspaces and Example~\ref{exa:c1nonc2} shows that 2D 
projections of sets from $\cC_3'$ may not lie in $\cC_2^0$. 
\par
Throughout this section let $K\subset\bR^m$, $m\in\bN$, be a convex body 
including the origin $0\in\bR^m$ in its interior. The {\em polar} of $K$ is
\begin{equation}\label{def:polar}
K^\circ:=\{u\in\bR^m\mid\langle u,x\rangle\leq 1, x\in K\}.
\end{equation}
It is well-known \cite{Rockafellar1970,Schneider2014} that $K^\circ$ is 
a convex body with the origin in its interior and that 
$(K^\circ)^\circ=K$ holds. 
\begin{defn}\label{def:touching}
Let $C\subset\bR^m$ be a convex subset and assume $u\in\bR^m$ is such 
that $x\mapsto\langle x,u\rangle$ has a maximum on $C$. Then the
the exposed face $F:=\argmax_{x\in C}\langle x,u\rangle$ is non-empty
and the {\em touching cone} of $C$ at $u$ is defined as the face of the 
normal cone $N_C(F)$ including $u$ in its relative interior 
\cite{Schneider2014}. The normal cone $N_C(\emptyset)=\bR^m$ is a 
touching cone by definition. The set of touching cones of $C$ is denoted 
by $\cT_C$, the set of faces of $C$ by $\cF_C$.
\end{defn}
\par
Partially ordered by inclusion, $\cT_C$ and $\cF_C$ are complete lattices 
where the infimum is the intersection \cite{Weis2012a}. Recall that every 
normal cone is a touching cone and every exposed face is a face. We shall 
use Theorem~7.4 of \cite{Weis2012a}:
\begin{fact}\label{fa:T=face-of-N}
Let $C\subset\bR^m$ be a convex subset. Then $\cT_C$ is the set of 
non-empty faces of normal cones of $C$. If $C$ has an interior point
and if $T_0\in\cT_C$ is not equal to $\bR^m$ then the interval
$[\{0\},T_0]_{\cT_C}$ equals 
$\{T\subset\bR^m\mid \mbox{$T$ is a non-empty face of $T_0$} \}$.
\end{fact}
\par
The class $\cC_m''$ is closed under projection.
\begin{rem}\label{rem:projection}
Fact~\ref{fa:T=face-of-N} shows that a convex subset $C\subset\bR^m$, 
$m\in\bN$, with interior point lies in $\cC_m''$ if and only if 
$\cN_C=\cT_C$. This can be used to show that $\cC_m''$ is closed 
under projection. More precisely, for $k=1,\ldots,m$ any 
$k$-dimensional image of any element of $\cC_m''$ under a linear 
map $\bR^m\to\bR^k$ belongs to $\cC_k''$. Indeed, the property that 
all non-empty faces of normal cones are normal cones is passed from 
$C$ to its linear images by Corollary~7.7 of \cite{Weis2012a}.
\end{rem}
\par
The lattices of $K$ and $K^\circ$ are related by a commutative diagram 
\cite{Weis2012a,Weis2012b}:
\begin{equation}\label{eq:iso-polar}
\begin{tikzpicture}
\node (a) at (-3,\diagh) {$\cF_K$};
\node at (-2,\diagh) {$\supset$};
\node (b) at (-1,\diagh) {$\cE_K$};
\node at (0,\diagh) {$\cong$};
\node (c) at (1,\diagh) {$\cN_K$};
\node at (2,\diagh) {$\subset$};
\node (d) at (3,\diagh) {$\cT_K$};
\node (e) at (-3,-\diagh) {$\cT_{K^\circ}$};
\node at (-2,-\diagh) {$\supset$};
\node (f) at (-1,-\diagh) {$\cN_{K^\circ}$};
\node at (0,-\diagh) {$\cong$};
\node (g) at (1,-\diagh) {$\cE_{K^\circ}$};
\node at (2,-\diagh) {$\subset$};
\node (h) at (3,-\diagh) {$\cF_{K^\circ}$};
\draw
(a) edge[->,font=\scriptsize,>=angle 90] node[left]{$\pos$} (e)
(b) edge[->,font=\scriptsize,>=angle 90] node[left]{$\pos$} (f)
(g) edge[->,font=\scriptsize,>=angle 90] node[right]{$\pos$} (c)
(h) edge[->,font=\scriptsize,>=angle 90] node[right]{$\pos$} (d);
\end{tikzpicture}
\end{equation}
The antitone lattice isomorphism (\ref{eq:exp-norm-iso}) appears  
once as $\cE_K\to\cN_K$ and once as $\cE_{K^\circ}\to\cN_{K^\circ}$. 
The positive hull operator $\pos$ defines isotone lattice isomorphisms 
$\cF_K\to\cT_{K^\circ}$, $\cE_K\to\cN_{K^\circ}$, 
$\cE_{K^\circ}\to\cN_K$, and $\cF_{K^\circ}\to\cT_K$ (we set 
$\pos(\emptyset)=\{0\}$). The map $\cE_K\to\cE_{K^\circ}$ which makes 
the diagram commute maps $F\in\cE_K$ to its {\em conjugate face}
$\{u\in K^\circ\mid\langle u,x\rangle=1\mbox{ for all }x\in F\}$.
\par
The diagram (\ref{eq:iso-polar}) allows Definition~\ref{def:main}
to be reformulated in terms of the polar. 
\begin{thm}\label{thm:def-main}~
\begin{enumerate}
\item
$K\in\cC_m^0$ $\iff$ every proper exposed face of $K^\circ$ 
contains an exposed point of $K^\circ$,
\item
$K\in\cC_m$ $\iff$ every proper exposed face $F$ of $K^\circ$ 
contains $\dim(F)+1$ affinely independent exposed points of $K^\circ$, 
\item
$K\in\cC_m'$ $\iff$ $K^\circ$ has no non-exposed points,
\item
$K\in\cC_m''$ $\iff$ $K^\circ$ has no non-exposed faces
$\iff$ $\cN_K=\cT_K$.
\end{enumerate} 
\end{thm}
{\em Proof:} 
Let $N$ be a proper normal cone of $K$. Then $N$ is a pointed closed cone. 
Let $F$ be the unique proper exposed face of $K^\circ$ such that 
$N=\pos(F)$, which exists by the isomorphism $\cE_{K^\circ}\to\cN_K$ of 
diagram (\ref{eq:iso-polar}). Then $\pos$ restricts to a lattice 
isomorphism $\cF_F\to\cF_N\setminus\{\emptyset\}$ as was observed in 
Lemma~3.4 of \cite{Weis2012a}. It is easy to see that $\pos$ restricts 
further to a lattice isomorphism $\cE_F\to\cE_N\setminus\{\emptyset\}$.
\par
The isomorphism $\cE_{K^\circ}\to\cN_K$ of diagram (\ref{eq:iso-polar}),
which identifies proper exposed faces $F$ of $K^\circ$ with proper normal 
cones $N$ of $K$ combined with the isomorphisms 
$\cF_F\to\cF_N\setminus\{\emptyset\}$ and 
$\cE_F\to\cE_N\setminus\{\emptyset\}$ proves (1), (2), and proves further
that $K\in\cC_m''$ holds if and only if every face of every proper exposed 
face of $K^\circ$ lies in $\cE_{K^\circ}$. 
\par
To prove (4) we assume the last statement to be true and show that any 
proper face $G$ of $K^\circ$ is exposed. The proper face $G$ lies in a 
proper exposed face $F$ of $K^\circ$, hence $G$ is a face of $F$ and by 
assumption $G$ is an exposed face of $K^\circ$. This shows 
$\cE_{K^\circ}=\cF_{K^\circ}$. The converse is clear because every face
of every face of $K^\circ$ is a face of $K^\circ$ \cite{Rockafellar1970}. 
The second condition $\cN_K=\cT_K$ of (4) follows from the isomorphisms 
$\cE_{K^\circ}\to\cN_K$ and $\cF_{K^\circ}\to\cT_K$ of diagram 
(\ref{eq:iso-polar}) or from Fact~\ref{fa:T=face-of-N} (already observed 
in Remark~\ref{rem:projection}). 
\par
The proof of (3) is analogous to the proof of the first condition of (4), 
now with $G$ an extreme point rather than a face.
\hspace*{\fill}$\square$\\
\par
We observe that Theorem~\ref{thm:def-main}(4) simplifies 
(\ref{eq:iso-polar}) to the following commutative diagram, valid 
for $K\in\cC_m''$.
\begin{equation}\label{eq:iso-polar-best}
\begin{tikzpicture}
\node (a) at (-3,\diagh) {$\cF_K$};
\node at (-2,\diagh) {$\supset$};
\node (b) at (-1,\diagh) {$\cE_K$};
\node at (0,\diagh) {$\cong$};
\node (c) at (1,\diagh) {$\cN_K$};
\node at (2,\diagh) {$=$};
\node (d) at (3,\diagh) {$\cT_K$};
\node (e) at (-3,-\diagh) {$\cT_{K^\circ}$};
\node at (-2,-\diagh) {$\supset$};
\node (f) at (-1,-\diagh) {$\cN_{K^\circ}$};
\node at (0,-\diagh) {$\cong$};
\node (g) at (1,-\diagh) {$\cE_{K^\circ}$};
\node at (2,-\diagh) {$=$};
\node (h) at (3,-\diagh) {$\cF_{K^\circ}$};
\draw
(a) edge[->,font=\scriptsize,>=angle 90] node[left]{$\pos$} (e)
(b) edge[->,font=\scriptsize,>=angle 90] node[left]{$\pos$} (f)
(g) edge[->,font=\scriptsize,>=angle 90] node[right]{$\pos$} (c)
(h) edge[->,font=\scriptsize,>=angle 90] node[right]{$\pos$} (d);
\end{tikzpicture}
\end{equation}
Notably, if $K\in\cC_m''$ then Minkowski's theorem shows that 
$\cE_{K^\circ}=\cF_{K^\circ}$ is atomistic, so $\cE_K$ is coatomistic 
(compare to Corollary~\ref{cor:add} for general convex sets). 
\par
Another consequence of the isomorphism $\cE_{K^\circ}\to\cN_K$ of 
diagram (\ref{eq:iso-polar}) and of Lemma~\ref{lem:char-strict}(2)
is the following.
\begin{cor}\label{cor:smooth-m2}
The convex body $K$ is smooth if and only if $K^\circ$ is 
strictly convex. 
\end{cor}
\par
In a sense, the next statement generalizes 
Corollary~\ref{cor:smooth-m2} from $\cC_m''$ to the class $\cC_m'$. 
The proof is inspired by Proposition~2.2 of \cite{Sanyal-etal2011}.
\begin{cor}\label{cor:coorbitope}
If the extreme points of $K^\circ$ lie on the boundary of a strictly 
convex set with interior points then $K\in\cC_m'$.
\end{cor}
{\em Proof:} 
Let $C$ be the strictly convex set containing the extreme points of 
$K^\circ$. Every boundary point of $C$ is an exposed point of $C$ by 
Lemma~\ref{lem:char-strict} so the extreme points of $K^\circ$ are 
exposed points of $C$ and {\em a fortiori} of $K^\circ$. This proves
that $K^\circ$ has no non-exposed points so
Theorem~\ref{thm:def-main}(3) shows $K\in\cC_m'$. 
\hspace*{\fill}$\square$\\
\par
We finish the discussion of Remark~\ref{rem:def}. 
\begin{defn}\label{def:minor}
Let $\cC_m'^*$ denote the class of proper convex subsets of $\bR^m$, $m\in\bN$, 
such that every exposed ray of every proper normal cone of $C$ lies in $\cN_C$.
\end{defn}
\par
Lemma~\ref{lem:max-number-rays} proves $\cC_m'^* \subset \cC_m$ and this implies 
\begin{equation}\label{eq:nesting}
\cC_m'' \subset \cC_m' \subset \cC_m'^* \subset \cC_m \subset \cC_m^0.
\end{equation}
For $m=3$ the first three inclusions are strict by the
Examples~\ref{exa:c1nonc2}, \ref{exa:AltCm1}, \ref{exa:cnonc1}. The last 
inclusion is strict already for $m=2$ by Example~\ref{exa:intro}.
\begin{exa}[Convex hull of ball and lens, $\cC_3''\subsetneq\cC_3'$]
\label{exa:c1nonc2}
We use Corollary~\ref{cor:coorbitope} and construct a convex body $K$
the extreme points of whose polar $K^\circ$ lie on a sphere. Consider 
the lens $C_{\rm ()}=D_-\cap D_+$ from Example~\ref{exa:intro}, where
\[
D_\pm:=\{(x,y)\in\bR^2\mid(x\pm\tfrac{3}{2})^2+y^2\leq(\tfrac{5}{2})^2\},
\]
and its embedding into $\bR^3$ defined by
\[
\widetilde{C_{\rm ()}}
:=C_{\rm ()}\oplus\{0\}
=\{(x,y,0)\in\bR^3\mid (x,y)\in C_{\rm ()}\}.
\]
Let $K:=\conv(\widetilde{C_{\rm ()}}\cup B)$ where $B$ is the 
closed Euclidean unit ball of $\bR^3$. Corollary~16.5.2 of 
\cite{Rockafellar1970} proves 
\begin{equation}\label{eq:c1nonc0}
K^\circ
=\conv(D_-^\circ\cup D_+^\circ)\oplus\bR\,\cap\, B
\end{equation}
which shows clearly that the extreme points of 
$K^\circ$ lie on the unit sphere of $\bR^3$. So 
Corollary~\ref{cor:coorbitope} proves $K\in\cC_3'$. 
\par
We use Theorem~5.23 and Algorithm~5.1 of 
\cite{RostalskiSturmfels2013} to compute the polar of $D_\pm$,
\[
D_\pm^\circ
=\{(x,y)\in\bR^2\mid \tfrac{1}{25} (8 x \mp 3)^2 + (2 y)^2 \leq 1\}.
\]
Clearly $\conv(D_-^\circ\cup D_+^\circ)$ has the non-exposed points 
$(\sigma\tfrac{3}{8},\tau\tfrac{1}{2})$ for signs $\sigma,\tau\in\{+,-\}$. 
The two ellipses $D_\pm^\circ$ lie in the unit disk, so the segments
\[
[(\sigma\tfrac{3}{8},\tau\tfrac{1}{2},-\tfrac{\sqrt{39}}{8}),
(\sigma\tfrac{3}{8},\tau\tfrac{1}{2},+\tfrac{\sqrt{39}}{8})],
\quad
\sigma,\tau\in\{+,-\}
\]
are non-exposed faces of $K^\circ$ by equation (\ref{eq:c1nonc0}). 
Now Theorem~\ref{thm:def-main}(4) shows $K\not\in\cC_3''$.
\par
The projection of $K$ onto the $x$-$y$-plane is $C_{\rm ()}$ which 
does not belong to $\cC_2^0$. Therefore the classes $\cC_m^0$, $\cC_m$, 
and $\cC_m'$ are not closed under linear maps.
\end{exa}
\begin{figure}
\includegraphics[height=2cm]{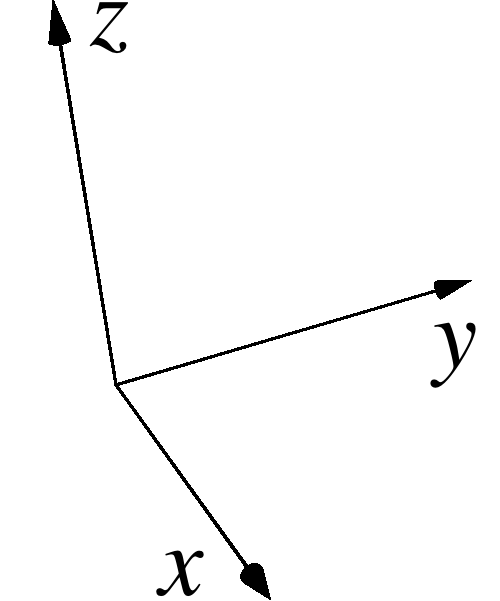}
a)\includegraphics[height=4cm]{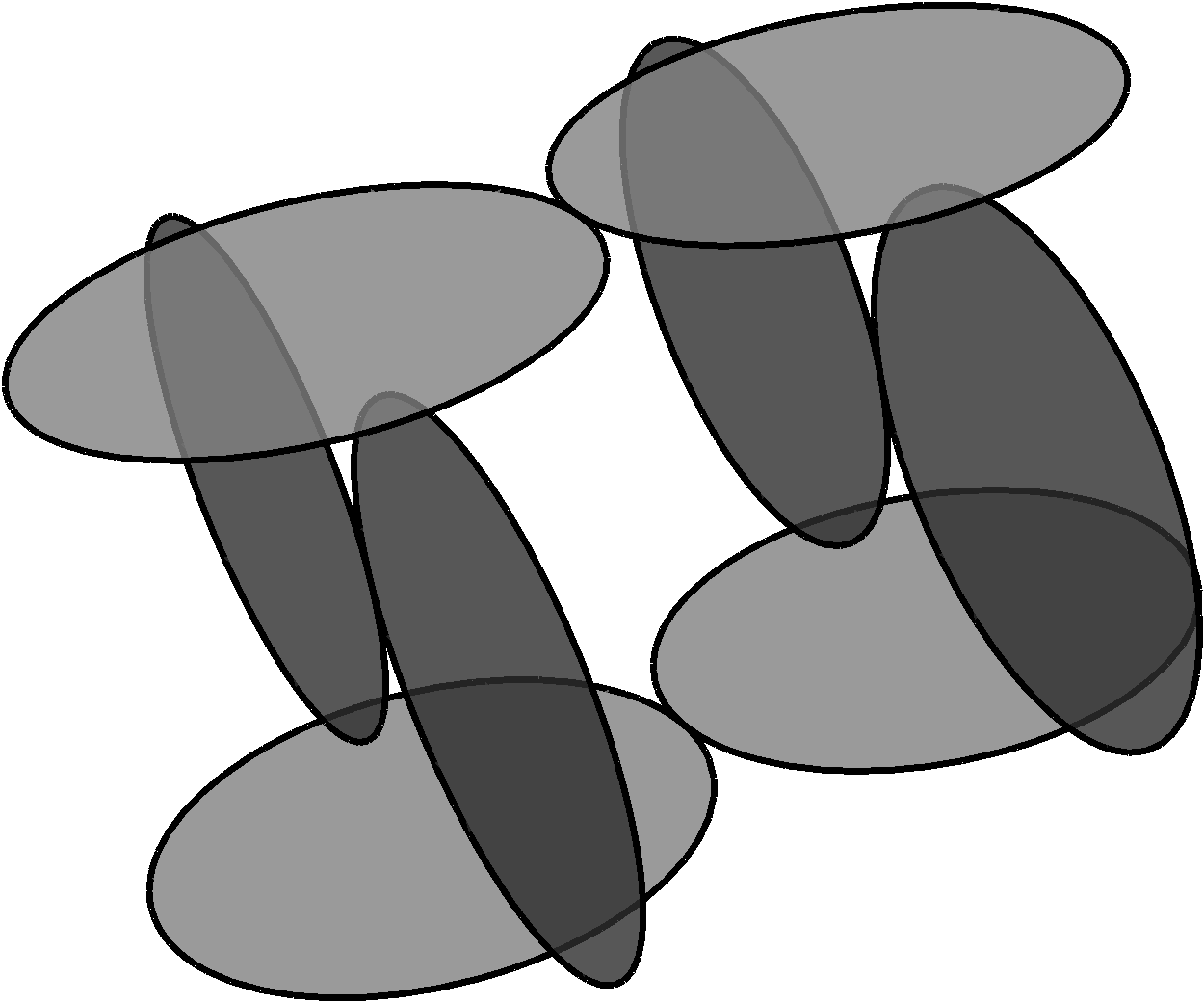}
b)\includegraphics[height=4cm]{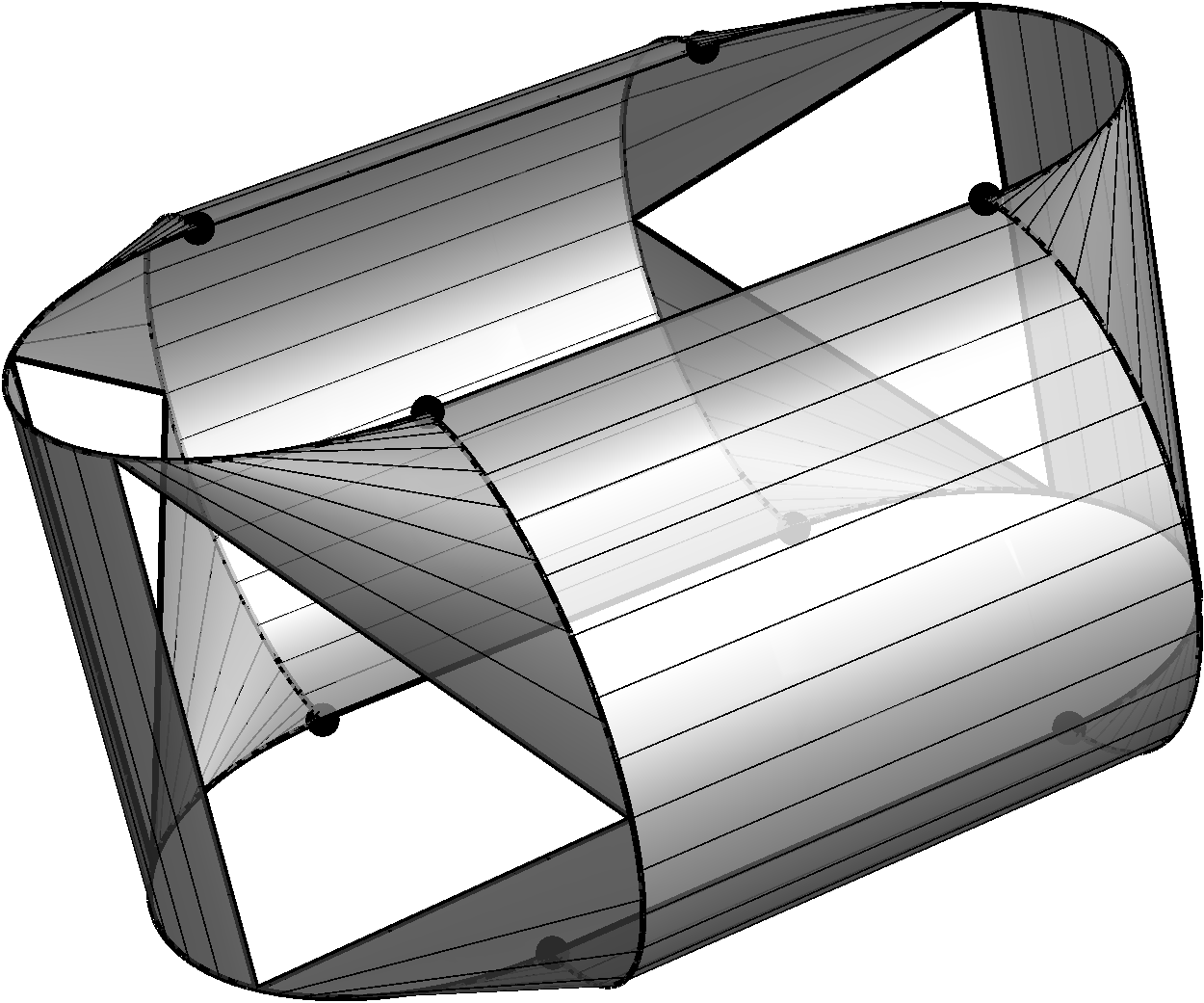}
\caption{\label{fig:AltCm1}% 
a) Unit disks with normal vectors in $z$-direction centered at 
$(0,\pm1,\pm1)$ (respectively, $y$-direction at $(\pm1,\pm1,0)$).
b) The convex hull of the eight disks from a) is sketched, only 
the boundary is depicted, two-dimensional faces are deleted, 
ruled surfaces are striped, and non-exposed points are marked at 
$(\pm1,\pm1,\pm1)$.}
\end{figure}
\begin{exa}[Polar of four stadia, $\cC_3' \subsetneq \cC_3'^*$]\label{exa:AltCm1}
Analogous to Theorem~\ref{thm:def-main}(3) one proves that $K\in\cC_m'^*$ 
holds if and only if every exposed point of every proper exposed face of 
$K^\circ$ lies in $\cE_{K^\circ}$. Therefore, if $L$ is a convex body with a 
non-exposed point, every exposed point of whose proper exposed faces is an 
exposed point of $L$ then the polar $K^\circ=L$ satisfies
$K\in\cC_3'^*\setminus\cC_3'$.
\par
Let $L$ be the convex set defined in Figure~\ref{fig:AltCm1}, which
is the convex hull of a cube with four half-cylinders attached, or the 
convex hull of four stadia. Two of the stadia lie in $\cE_L$
and include eight non-exposed points of $L$. Their four straight sides 
are non-exposed faces of $L$. Therefore the eight non-exposed points 
of $L$ are no exposed points of any proper exposed face of $L$. This
completes the example.
\par
It has nothing to do with the inclusion $\cC_3' \subsetneq \cC_3'^*$,
but is nevertheless interesting that $L$ is the convex support 
(\ref{def:convex-support}) of three
$16$-by-$16$ matrices and therefore the results of the 
Sections~\ref{sec:smoothness} and~\ref{sec:proof-main} apply to $L$. 
Indeed, the unit disk with normal vector in 
$z$-direction centered at $(0,1,1)$ is the convex support of 
\[
F_1=\left(\begin{smallmatrix}0 & 1\\1 & 0\end{smallmatrix}\right),
\qquad
F_2=\left(\begin{smallmatrix}1 & -\ii\\\ii & 1\end{smallmatrix}\right),
\qquad
F_3=\left(\begin{smallmatrix}1 & 0\\0 & 1\end{smallmatrix}\right).
\]
Similarly the other seven disks of Figure~\ref{fig:AltCm1}a) are 
convex support sets and their convex hull $L$ is the convex support 
of the direct sum of those matrices.
\par
The lattice of exposed faces $\cE_L$ has six coatoms of dimension two 
(two stadia and four triangles). The reminder of the boundary of $L$ is 
covered by four cylindrical and eight curved ruled surfaces. The straight 
segments on the curved ruled surface of the positive octant have 
end-points $(\cos(\varphi), 1+\sin(\varphi),1)$ and 
$(1+\sin(\psi),1,\cos(\psi))$, where
\[
\cos(\psi)=\tfrac{2\cos(\varphi)-1}{2-2\cos(\varphi)+\cos(2\varphi)},
\qquad
\sin(\psi)=\tfrac{4\cos(\varphi)\sin(\tfrac{\varphi}{2})^2}{%
2-2\cos(\varphi)+\cos(2\varphi)},
\qquad
\varphi\in(0,\tfrac{\pi}{3}).
\]
All boundary points of $L$ are smooth except for the exposed points. 
They cover a half-circle of each of the eight disks used to define 
$L$ in Figure~\ref{fig:AltCm1}. Having a two-dimensional normal cone,
each of them is the intersection of a unique pair of coatoms of $\cE_L$
(segment, triangle, or stadium) by 
Corollary~\ref{cor:max-number-coatoms}(2).
\end{exa}
\begin{exa}[Polar of truncated torus, $\cC_3'^*\subsetneq\cC_3$]\label{exa:cnonc1}
Consider the convex hull $K$ of the torus
$\{(x,y,z)\in\bR^3\mid p(x,y,z)=0\}$ defined by 
\[
p:=(x^2 + y^2 + z^2)^2-4 (x^2 + y^2)
\]
which is a torus with rotation symmetry about the $z$-axis and
self-intersection at the origin $0\in\bR^3$. The radius from $0$
to the center of the torus tube and the radius of the tube are 
both equal one. The convex body $K$ is smooth, so $K\in\cC_3''$,
and has the following proper exposed faces:
\[
D_\pm:=\{(x,y,\pm1)\in\bR^3\mid x^2+y^2\leq 1\}
\qquad\mbox{(two disks)}
\]
and exposed points $(x,y,z)$ for $x,y,z\in\bR$ such that $p(x,y,z)=0$
and $x^2+y^2>1$. The non-exposed points of $K$ are 
\[
\xi_\pm(\varphi):=(\cos(\varphi),\sin(\varphi),\pm1),
\quad
\varphi\in[0,2\pi).
\]
By intersecting $K$ with three half-spaces 
($\varphi_0=0,\varphi_1=\tfrac{2}{3}\pi,\varphi_2=\tfrac{4}{3}\pi$)
we define
\[
L:=\{(x,y,z)\in K\mid x\cos(\varphi_i)+y\sin(\varphi_i)\leq 1,i=0,1,2\}.
\]
The convex body $L$ has five 2D exposed faces, namely $D_\pm$ and 
the exposed faces in the boundaries of the three half-spaces. The 
latter three exposed faces of $L$ contain no non-exposed points of 
$L$. The disk $D_\pm$ contains three exposed points 
$\xi_\pm(\varphi_1)$, $\xi_\pm(\varphi_2)$, $\xi_\pm(\varphi_3)$
of $L$ while the remaining exposed points of $D_\pm$ are 
non-exposed points of $L$. The convex body  $L$ has no 1D exposed 
faces. This discussion of exposed points of $L$ and 
Theorem~\ref{thm:def-main}(2),(3) show 
$L^\circ\in\cC_3\setminus\cC_3'$.
\end{exa}
%
%%%%%%%%%%%%%%%%%%%%%%%%%%%%%%%%%%%%%%%%%%%%%%%%%%%%%%%%%%%%%%%%%%%%%%%%%%%%
%%%%%%%%%%%%%%%%%%%%%%%%%%%%%%%%%%%%%%%%%%%%%%%%%%%%%%%%%%%%%%%%%%%%%%%%%%%%
%%%%%%%%%%%%%%%%%%%%%%%%%%%%%%%%%%%%%%%%%%%%%%%%%%%%%%%%%%%%%%%%%%%%%%%%%%%%
%%%%%%%%%%%%%%%%%%%%%%%%%%%%%%%%%%%%%%%%%%%%%%%%%%%%%%%%%%%%%%%%%%%%%%%%%%%%
%%%%%%%%%%%%%%%%%%%%%%%%%%%%%%%%%%%%%%%%%%%%%%%%%%%%%%%%%%%%%%%%%%%%%%%%%%%%
%
\section{State spaces of operator systems}
\label{sec:state-spaces}
\par
Theorem~\ref{thm:state-space} proves
that the state space $\cM(S)$ of an operator system $S$ of $M_n$ 
is the projection of the state space $\cM_n$ of the matrix algebra 
$M_n$ onto the real space $S\her\subset S$ of hermitian matrices. 
We point out that
\begin{itemize}
\item $\cM(S)\in\cC_m''$ for $m=\dim_\bC(S)-1$
(remove one codimension of $\cM(S)$ from $S\her$),
\item
there is a lattice isomorphism from the ground state projections of 
$S\her$ to the exposed faces of $\cM(S)$,
\item
the ground state projections of $S\her$ form a coatomistic lattice.
\end{itemize}
Some of these properties may be well-known, but we are unaware of proofs 
in the literature. We finish the article with representations of $\cM(S)$
as a convex support set and convex hull of a joint numerical range and 
with a discussion of convex support sets of 3-by-3 matrices. 
\par
Let $S\subset M_n$ be an operator system. It is well-known that every 
complex linear functional $f:S\to\bC$ has the form $f=f_a$ for a unique 
$a\in S$ where $f_a(s)=\langle a,s\rangle$ for all $s\in S$.
Here $\langle a,b\rangle=\tr(a^*b)$ is the Hilbert-Schmidt inner 
product of $a,b\in S$, whose real part defines a Euclidean scalar 
product on $S$ where $S=S\her+\ii S\her$ is an orthogonal direct 
sum. The functional $f_a$ has real values on $S\her$ if and only if 
$a\in S\her$ so the Hilbert-Schmidt inner product restricts to a
Euclidean scalar product on $S\her$.
\begin{defn}
A {\em state} on an operator system $S\subset M_n$ is a complex linear 
functional $f:S\to\bC$ which is {\em positive}, that is $f(s)\geq 0$ 
holds for all positive-semidefinite $s\in S$, and {\em normalized}, 
that is $f(\id)=1$. Let
\[
\cM(S):=\{s\in S\mid\mbox{$f_s$ is a state on $S$}\}.
\]
In abuse of notation we call elements of $\cM(S)$ {\em states} on $S$ 
and say that $\cM(S)$ is the {\em state space} of $S$.
\end{defn}
\par
Let $K\subset\bR^m$ be a convex body including the origin $0\in\bR^m$ 
in its interior. The {\em dual} of $K$ is defined as the point 
reflection of the polar $K^\circ$, defined in (\ref{def:polar}),
\[
K^*
:=-K^\circ
=\{u\in\bR^m\mid 1+ \langle u,x\rangle\geq 0, x\in K\}.
\]
It follows from the properties of $K^\circ$ that $K^*$ is a convex body 
with the origin in its interior and that $(K^*)^*=K$ holds. Similarly, 
the {\em dual} of a convex cone $K\subset\bR^m$ is 
\[
K^*=\{u\in\bR^m\mid \langle u,x\rangle\geq 0, x\in K\}
\]
and $K$ is {\em self-dual}, if $K=K^*$ holds. For every real subspace 
$X\subset(M_n)\her$ the mapping $\pi_X:(M_n)\her\to X$ denotes the 
orthogonal projection from $(M_n)\her$ onto $X$.
\par
The following proof uses a property of self-dual convex cones:
Projections of certain cone bases are dual to bounded affine sections 
through their interior.
\begin{thm}\label{thm:state-space}
If $S\subset M_n$ is an operator system then
$\cM(S)=\pi_{S\her}(\cM_n)$.
\end{thm}
{\em Proof:} 
It is well-known, see for instance \cite{BermanBen-Israel1973}, that the 
cone $K$ of positive semi-definite matrices of $M_n$ is self-dual with 
respect to $(M_n)\her$. Without loss of generality, let $\dim_\bC(S)>1$. 
Then the space $S_0\subset S\her$ of traceless hermitian matrices of $S$ 
has dimension $\dim_\bR(S_0)=\dim_\bC(S)-1>0$. The affine section
\[
\bA:=K\cap(\tfrac{1}{n}\id+S_0)=\cM_n\cap(\tfrac{1}{n}\id+S_0)
\]
of $K$ is bounded and contains the trace state $\tfrac{1}{n}\id$, 
which is an interior point with respect to $(M_n)\her$ of $K$. 
Since $\cM_n=(\tfrac{1}{n}\id+\id^\perp)\cap K$ holds, where $\id^\perp$ 
denotes orthogonal complement, Theorem~2.16 of \cite{Weis2011} proves 
(duals with respect to $S_0$)
\begin{enumerate}[label=\emph{\roman*)}, widest=iii]
\item $\bA-\tfrac{1}{n}\id=\tfrac{1}{n}\pi_{S_0}(\cM_n)^*$
\item $\pi_{S_0}(\cM_n)=\tfrac{1}{n}(\bA-\tfrac{1}{n}\id)^*$.
\end{enumerate}
On the other hand, an easy computation shows
\begin{align*}
\cM(S) &= \{ s\in S\her \mid 
\langle s,\id\rangle=1,
\langle s,a\rangle\geq 0 \, \forall a\in\bA\}\\
 &= \tfrac{1}{n}\id + \{ s_0\in S_0 \mid 
\langle \tfrac{1}{n}\id + s_0,a\rangle\geq 0 \, \forall a\in\bA\}\\
 &= \tfrac{1}{n}\id + \{ s_0\in S_0 \mid 
\tfrac{1}{n} + 
\langle s_0,a_0\rangle\geq 0 \, \forall a_0\in\bA-\tfrac{1}{n}\id\},
\end{align*}
that is
\begin{enumerate}[resume*]
\item
$\cM(S) = \tfrac{1}{n}(\id + (\bA-\tfrac{1}{n}\id)^*)$.
\end{enumerate}
Substituting {\em ii)} into {\em iii)} gives
$\cM(S)=\tfrac{1}{n}\id+\pi_{S_0}(\cM_n)$ which proves the claim.
\hspace*{\fill}$\square$\\
\par
State spaces lie in $\cC_m''$. More precisely the following holds,
if we identify the space of traceless hermitian matrices of a 
$d$-dimensional operator system with $\bR^{d-1}$.
\begin{cor}\label{cor:state-spaces-C2}
If $S\subset M_n$ is an operator system of dimension $d:=\dim_\bC(S)>1$ 
then $\cM(S)-\tfrac{1}{n}\id\in\cC_{d-1}''$. 
\end{cor}
{\em Proof:} 
One can see from the proof of Theorem~\ref{thm:state-space} 
that zero is an interior point with respect to $S_0$ of 
$\cM(S)-\tfrac{1}{n}\id=\pi_{S_0}(\cM_n)$. 
Further {\em i)} of Theorem~\ref{thm:state-space}  shows that the 
polar of $\pi_{S_0}(\cM_n)$ is an affine section of $\cM_n$. 
Since $\cM_n$ has no non-exposed faces 
\cite{RamanaGoldman1995,AlfsenShultz2001} 
it follows that this affine section has no non-exposed faces either. 
Then Theorem~\ref{thm:def-main}(4) proves the claim.
\hspace*{\fill}$\square$\\
\par
Analogues of Corollary~\ref{cor:state-spaces-C2} were proved in Fact~3.2 
of \cite{Szymanski-etal2016} for convex support sets, defined below. One 
of the proofs uses that $\cM_n-\tfrac{1}{n}\id$ lies in $\cC_{n^2-1}''$ 
so the projection $\pi_{S_0}(\cM_n)$ lies in $\cC_{d-1}''$ by 
Remark~\ref{rem:projection}. 
\par
We provide details about the lattice isomorphism from ground state 
projections of hermitian operators of an operator system to the 
exposed faces of the state space. 
\begin{rem}\label{rem:iso-exp-ground}
Let us begin with the full operator system $S=M_n$ where the 
isomorphism is well-known, see Section~6 of \cite{BarkerCarlson1975}, 
Chapter~3 of \cite{AlfsenShultz2001}, and the references therein.
We denote the set of projections $p=p^*=p^2\in M_n$ by $\cP$. Endowed 
with the partial ordering on $(M_n)\her$ defined by
\begin{equation}\label{eq:sd-order}
a\preceq b 
\quad\iff\quad
b\succeq a
\quad\iff\quad
\mbox{$b-a$ is positive semi-definite}
\end{equation}
the set $\cP$ is a complete lattice. One can represent $\cP$ in terms 
of the images of the projections, ordered by inclusion, see for example 
Chapter 2 of \cite{AlfsenShultz2001}. The infimum in this lattice of 
subspaces is the intersection. Let $u\in(M_n)\her$ be a hermitian matrix 
and
\[
u=\sum_\lambda \lambda p_{u,\lambda},
\quad
\id=\sum_\lambda p_{u,\lambda},
\quad
p_{u,\lambda}=p_{u,\lambda}^*=p_{u,\lambda}^2\in M_n
\]
its spectral decomposition ($\lambda$ extends over the eigenvalues of $u$). 
The {\em ground state projection} of $u$ is the spectral projection 
$p_-(u):=\min_\lambda p_{u,\lambda_-(u)}$ of the least eigenvalue 
$\lambda_-(u)$ of $u$.
The exposed face of $\cM_n$ with inward pointing normal vector 
$u\in(M_n)\her$
\[
F_{\cM_n}(u):=\argmin\{\langle\rho,u\rangle\mid\rho\in\cM_n\}
\]
is well-known to be expressible in the form
\[
F_{\cM_n}(u)=\{\rho\in\cM_n\mid {\rm Image}(\rho)\subset{\rm Image}(p_-(u))\}
\]
in terms of the ground state projection of $u$, and $\cP\to\cE_{\cM_n}$, 
$p\mapsto F_{\cM_n}(-p)$, is a lattice isomorphism to the exposed faces of 
$\cM_n$.
\par
Consider now the state space $\cM(S)$ of an operator system 
$S\subset M_n$. It is easy to see that for $s\in S\her$ the exposed face
\[
F_{\cM(S)}(s):=\argmin\{\langle\rho,s\rangle\mid\rho\in\cM(S)\}
\]
lifts to the corresponding exposed face of $\cM_n$, that is 
$F_{\cM_n}(s)=\pi_{S\her}|_{\cM_n}^{-1}(F_{\cM(S)}(s))$. Indeed, the 
map
\begin{equation}\label{eq:iso-lift}
\cE_{\cM(S)}\to
\cE_{\cM(S)}':=\{\pi_{S\her}|_{\cM_n}^{-1}(F)\mid F\in\cE_{\cM(S)}\}
\end{equation}
is a lattice isomorphism whose range $\cE_{\cM(S)}'$ is ordered by 
inclusion. Moreover, the infimum of $\cE_{\cM(S)}'$ is just the 
restriction of the infimum of $\cE_{\cM_n}$ to $\cE_{\cM(S)}'$, which 
is the intersection (see Proposition 5.6 of \cite{Weis2012a}). Let
\[
\cP(S):=\{p_-(s)\mid s\in S\her\}\cup\{0\}
\]
denote the set of ground state projections of $S\her$ endowed with the 
partial ordering (\ref{eq:sd-order}). Notice that $p_-(s)$ may not lie 
in $S$ for some $s\in S\her$.
Since (\ref{eq:iso-lift}) is a lattice isomorphism,
the lattice isomorphism $\cP\to\cE_{\cM_n}$ restricts to a lattice 
isomorphism $\cP(S)\to\cE_{\cM(S)}'$ and combines with (\ref{eq:iso-lift})
to a lattice isomorphism
\begin{equation}\label{eq:iso-pro-exp}
\cP(S)\to\cE_{\cM(S)},
\quad 
p\mapsto \pi_{S_h}(F_{\cM_n}(-p)).
\end{equation}
The isomorphism (\ref{eq:iso-pro-exp}) is described in Section~3.1 of 
\cite{Weis2011}. We stress that the infimum of $\cP(S)$ is the restriction
of the infimum of $\cP$ because the the infimum of $\cE_{\cM(S)}'$ is the 
restriction of the infimum of $\cE_{\cM_n}$.
\end{rem}
\par
If projections are represented in terms of their images then the 
infimum of $\cE_{\cM(S)}$ and of $\cP(S)$ are both given by the 
intersection. This is especially nice in the following corollary.
\begin{cor}\label{cor:coatomistic}
If $S\subset M_n$ is an operator system then the lattice of ground 
state projections $\cP(S)$ is coatomistic.
\end{cor}
{\em Proof:} 
This follows from Corollary~\ref{cor:state-spaces-C2}, 
Lemma~\ref{lem:coatom-intersection}, 
and Corollary~\ref{cor:add}(1).
\hspace*{\fill}$\square$\\
\par
We now recall two coordinate representations of state spaces 
(and show that polytopes belong to them). 
We define the {\em convex support} of $F_i\in(M_n)\her$, $i=1,\ldots,k$, 
$k\in\bN$, by
\begin{equation}\label{def:convex-support}
\cs(F_1,\ldots,F_k):=\{\langle\rho, F_i\rangle_{i=1}^k\mid\rho\in\cM_n\}
\end{equation}
and the {\em joint numerical range}
(using the inner product of $\bC^n$) by
\[
W(F_1,\ldots,F_k):=\{\langle x,F_i(x)\rangle_{i=1}^k\mid \langle x,x\rangle=1,
x\in\bC^n\}.
\]
\begin{rem}[Convex support sets]\label{rem:convex-support}~
\begin{enumerate}
\item
{\em State spaces and convex support sets}.
Let $S\subset M_n$ be an operator system and let $F_i\in(M_n)\her$ for
$i=1,\ldots,k$, such that $\{\id,F_1,\ldots,F_k\}$ spans
$S\her$. Then it is easy to see that the linear map 
\[
\alpha:S\her\to\bR^{k+1},
s\mapsto(\langle s,\id\rangle,\langle s,F_1\rangle,\ldots,\langle s,F_k\rangle)
\]
restricts to a bijection $S\her\to\alpha(S\her)$. Hence
$\cM(S)\cong\cs(\id,F_1,\ldots,F_k)$ and 
$\cM(S)\cong\cs(F_1,\ldots,F_k)$ are isomorphic, see Remark~1.1(1) of
\cite{Weis2011} for a proof.
\item
{\em Convex support sets and joint numerical ranges}.
If $F_i\in(M_n)\her$, $i=1,\ldots,k$, then 
\[
\cs(F_1,\ldots,F_k)=\conv\left(W(F_1,\ldots,F_k)\right).
\] 
For a proof see \cite{Mueller2010} or Section~2 of 
\cite{Szymanski-etal2016}. For $k=2$ the joint numerical range 
$W(F_1,F_2)$ is the numerical range of $F_1+\ii F_2$ which is convex 
for all $n\in\bN$ by the Toeplitz-Hausdorff theorem. The joint 
numerical range is convex for $k=3$ and $n\geq 3$, but fails to be 
convex in general for $k\geq 4$ 
\cite{Au-YeungPoon1979,LiPoon2000,Gutkin-etal2004}.
\item
{\em Polytopes}. Consider a collection of $n$ points $p_j\in\bR^m$, 
$n,m\in\bN$. Define an $m$-by-$n$-matrix $M$ whose $j$'s column is 
$p_j$, $j=1,\ldots,n$, and define diagonal matrices $F_i\in M_n$ 
where $\diag(F_i)$ is the $i$'s row of $M$, $i=1,\ldots,m$. Then 
the convex support $\cs(F_1,\ldots,F_m)$ is the convex hull of 
$\{p_1,\ldots,p_n\}\subset\bR^m$.
\end{enumerate}
\end{rem}
\par
We discuss convex support sets of $3$-by-$3$ matrices.
\begin{exa}\label{exa:JNR3x3}
Consider the convex support $\cs(F)$ of $k\in\bN$ hermitian 
$3$-by-$3$ matrices $F=(F_1,\ldots,F_k)$ and assume without loss 
of generality that $\cs(F)$ has an interior point.
Then $\cs(F)\in\cC_k''$ holds by Remark~\ref{rem:convex-support}(1) 
and Corollary~\ref{cor:state-spaces-C2}.
\par
The subset of non-smooth points $\partial_{\rm n}\cs(F)$ of the boundary 
$\partial\cs(F)$ consists of the non-smooth exposed points. Indeed, 
The form of a proper exposed face $G$ of $\cs(F)$ is very restricted 
for $3$-by-$3$ matrices. It may be a singleton (exposed point), either 
smooth or not. Otherwise $G$ is a segment, a filled ellipse, or filled 
ellipsoid \cite{Szymanski-etal2016}. We see that proper exposed faces of 
$\cs(F)$ have no segments on their relative boundary. Therefore, if $G$ 
is not a singleton then it is necessarily a coatom of $\cE_{\cs(F)}$
and is therefore smooth (has a unique unit normal vector) by
Theorem~\ref{thm:one-ray}. For the same reason, every 
non-exposed face of $\cs(F)$ is a singleton. Moreover, every non-exposed 
point is smooth because it has the same normal cone as the smallest 
exposed face in which it is contained (see Lemma~4.6 of \cite{Weis2012a}). 
This proves the claim. Equivalently, $\partial_{\rm n}\cs(F)$ is the set 
of intersection points of pairs of mutually distinct coatoms of 
$\cE_{\cs(F)}$ by Corollary~\ref{cor:main}.
\par
A notable property of $3$-by-$3$ matrices is that each two coatoms of 
$\cE_{\cs(F)}$ which are both no singletons must intersect. This follows 
from the analogue property of $k=2$ by projecting $\cs(F)$ onto the span 
of the normal vectors of these coatoms  (for a proof when $k=3$ see 
Lemma~5.1 of \cite{Szymanski-etal2016}, $k>3$ is analogous). 
The case $k=2$ is solved by the classification of the numerical range of 
$3$-by-$3$ matrices \cite{Kippenhahn1951,Keeler-etal1997}. So, 
$\cE_{\cs(F)}$ has only one cluster of coatoms of positive dimension, in 
the sense of Corollary~\ref{cor:clustering}. 
\par
The classification \cite{Szymanski-etal2016} proves that 
$\partial_{\rm n}\cs(F)$ is closed for $k=3$. If $\cs(F)$ has a corner 
then $\cs(F)$ is either the convex hull of an ellipsoid and a point 
outside the ellipsoid or the convex hull of an ellipse and a point 
outside the affine hull of the ellipse. In the first case 
$\partial_{\rm n}\cs(F)$ is a singleton, in the second case 
$\partial_{\rm n}\cs(F)$ is the union of an ellipse and a singleton. 
If $\cs(F)$ has no corner then the cluster of $\cE_{\cs(F)}$ with 
coatoms of positive dimension contains $s$ segments and $e$ ellipses 
such that 
\[
(s,e)\in\{%
\begin{smallmatrix}(0,0)\end{smallmatrix},
\begin{smallmatrix}(0,1)\end{smallmatrix},
\begin{smallmatrix}(0,2)\end{smallmatrix},
\begin{smallmatrix}(0,3)\end{smallmatrix},
\begin{smallmatrix}(0,4)\end{smallmatrix},
\begin{smallmatrix}(1,0)\end{smallmatrix},
\begin{smallmatrix}(1,1)\end{smallmatrix},
\begin{smallmatrix}(1,2)\end{smallmatrix}\}.
\]
The coatoms of this cluster intersect in pairs but not in triples so 
$\partial_{\rm n}\cs(F)$ has cardinality ${s+e\choose 2}\in\{0,1,3,6\}$. 
Since $\partial_{\rm n}\cs(F)$ is closed it follows that its complement
$\partial\cs(F)\setminus\partial_{\rm n}\cs(F)$ is a $C^1$-submanifold 
of $\bR^3$ (Theorem 2.2.4 of \cite{Schneider2014} can be used locally
to prove this). 
\end{exa}
%
%
%%%%%%%%%%%%%%%%%%%%%%%%%%%%%%%%%%%%%%%%%%%%%%%%%%%%%%%%%%%%%%%%%%%%%%%%%%%%
%%%%%%%%%%%%%%%%%%%%%%%%%%%%%%%%%%%%%%%%%%%%%%%%%%%%%%%%%%%%%%%%%%%%%%%%%%%%
%%%%%%%%%%%%%%%%%%%%%%%%%%%%%%%%%%%%%%%%%%%%%%%%%%%%%%%%%%%%%%%%%%%%%%%%%%%%
%%%%%%%%%%%%%%%%%%%%%%%%%%%%%%%%%%%%%%%%%%%%%%%%%%%%%%%%%%%%%%%%%%%%%%%%%%%%
%%%%%%%%%%%%%%%%%%%%%%%%%%%%%%%%%%%%%%%%%%%%%%%%%%%%%%%%%%%%%%%%%%%%%%%%%%%%
%
\vspace{\baselineskip}
{\par\noindent\footnotesize
{\em Acknowledgements.}
The author thanks Eduardo Garibaldi, Fernando Torres, and especially
Marcelo Terra Cunha for discussions. He is grateful for Arleta Szko{\l}a's 
comments on an earlier manuscript. This work is supported by a PNPD/CAPES 
scholarship of the Brazilian Ministry of Education.}
%
%%%%%%%%%%%%%%%%%%%%%%%%%%%%%%%%%%%%%%%%%%%%%%%%%%%%%%%%%%%%%%%%%%%%%%%%%%%%
%%%%%%%%%%%%%%%%%%%%%%%%%%%%%%%%%%%%%%%%%%%%%%%%%%%%%%%%%%%%%%%%%%%%%%%%%%%%
%%%%%%%%%%%%%%%%%%%%%%%%%%%%%%%%%%%%%%%%%%%%%%%%%%%%%%%%%%%%%%%%%%%%%%%%%%%%
%%%%%%%%%%%%%%%%%%%%%%%%%%%%%%%%%%%%%%%%%%%%%%%%%%%%%%%%%%%%%%%%%%%%%%%%%%%%
%%%%%%%%%%%%%%%%%%%%%%%%%%%%%%%%%%%%%%%%%%%%%%%%%%%%%%%%%%%%%%%%%%%%%%%%%%%%
%
\bibliographystyle{plain}

\begin{thebibliography}{10}
%
\bibitem{AlfsenShultz2001} E.\,M.~Alfsen and F.\,W.~Shultz (2001)
{\em State Spaces of Operator Algebras:
Basic Theory, Orientations, and C*-Products},
Boston: Birkh\"auser
%
\bibitem{Arrachea-etal1992} L.~Arrachea, N.~Canosa, A.~Plastino, 
M.~Portesi, and R.~Rossignoli (1992)
{\em Maximum-entropy approach to critical phenomena in ground states 
of finite systems},
Phys Rev A {\bf 45} 7104--7110
%
\bibitem{Au-YeungPoon1979} Y.\,H.~Au-Yeung and Y.\,T.~Poon (1979)
{\em A remark on the convexity and positive definiteness
concerning Hermitian matrices},
Southeast Asian Bull Math {\bf 3} 85--92
%
\bibitem{BarkerCarlson1975} G.\,P.~Barker and D.~Carlson (1975)
{\em Cones of diagonally dominant matrices},
Pac J Math {\bf 57} 15--32
%
\bibitem{Barker1978} G.\,P.~Barker (1978)
{\em Faces and duality in convex cones},
Linear Multilinear A {\bf 6} 161--169
%
\bibitem{Barndorff-Nielsen1978} O.~Barndorff-Nielsen (1978)
{\em Information and Exponential Families in Statistical Theory},
Chichester: John Wiley \& Sons Ltd 
%
\bibitem{BengtssonZyczkowski2006} I.~Bengtsson and K.~\.Zyczkowski (2006)
{\em Geometry of Quantum States}, Cambridge: Cambridge University Press
%
\bibitem{BermanBen-Israel1973} A.~Berman and A.~Ben-Israel (1973)
{\em Linear equations over cones with interior: A solvability theorem 
with applications to matrix theory},
Linear Algebra Appl {\bf 7} 139--149
%
\bibitem{Birkhoff1967} G.~Birkhoff (1967)
{\em Lattice Theory}, 3rd ed., Providence, R.I.: AMS
%
\bibitem{Chen-etal2012a} J.~Chen, Z.~Ji, D.~Kribs, Z.~Wei, and B.~Zeng (2012)
{\em Ground-state spaces of frustration-free Hamiltonians},
J Math Phys {\bf 53} 102201
%
\bibitem{Chen-etal2015} J.~Chen, Z.~Ji, C.-K.~Li, Y.-T.~Poon, Y.~Shen, N.~Yu, 
B.~Zeng, and D.~Zhou (2015)
{\em Discontinuity of maximum entropy inference and quantum phase transitions},
New J Phys {\bf 17} 083019
%
\bibitem{Chen-etal2012b} J.~Chen, Z.~Ji, M.\,B.~Ruskai, B.~Zeng, 
and D.-L.~Zhou (2012) 
{\em Comment on some results of Erdahl and the convex structure of 
reduced density matrices},
J Math Phys {\bf 53} 072203
%
\bibitem{Chen-etal2012c} J.~Chen, Z.~Ji, B.~Zeng, and D.\,L.~Zhou (2012)
{\em From ground states to local Hamiltonians}, Phys Rev A {\bf 86} 022339
%
\bibitem{Cheung-etal2011} W.-S.~Cheung, X.~Liu, and T.-Y.~Tam (2011)
{\em Multiplicities, boundary points, and joint numerical ranges},
Oper Matrices {\bf 1} 41--52
%
\bibitem{ChienNakazato2010} M.-T.~Chien and H.~Nakazato (2010)
{\em Joint numerical range and its generating hypersurface},
Linear Algebra Appl {\bf 432} 173--179
%
\bibitem{Erdahl1972} R.\,M.~Erdahl (1972)
{\em The convex structure of the set of N-representable reduced 2-matrices},
J Math Phys {\bf 13} 1608--1621
%
\bibitem{Gruenbaum2003} B.~Gr\"unbaum (2003)
{\em Convex Polytopes}, 2nd ed., New York: Springer
%
\bibitem{Gutkin-etal2004} E.~Gutkin, E.\,A.~Jonckheere, and M.~Karow (2004)
{\em Convexity of the joint numerical range: topological and 
differential geometric viewpoints},
Linear Algebra Appl {\bf 376} 143--171
%
\bibitem{GutkinZyczkowski2013} E.~Gutkin and K.~\.Zyczkowski (2013)
{\em Joint numerical ranges, quantum maps, and joint numerical shadows},
Linear Algebra Appl {\bf 438} 2394--2404
%
\bibitem{Heinosaari-etal2013} T.~Heinosaari, L.~Mazzarella, 
and M.\,M.~Wolf (2013)
{\em Quantum tomography under prior information},
Commun Math Phys {\bf 318} 355--374
%
\bibitem{Holevo2011} A.~Holevo (2011) 
{\em Probabilistic and Statistical Aspects of Quantum Theory},
2nd edition, Pisa: Edizioni della Normale
%
\bibitem{Jaynes1957} E.~Jaynes (1957)
{\em Information theory and statistical mechanics. II},
Phys Rev {\bf 108} 171--190
%
\bibitem{Keeler-etal1997} D.\,S.~Keeler, L.~Rodman, and
I.\,M.~Spitkovsky (1997)
{\em The numerical range of $3\times 3$ matrices},
Lin Alg Appl {\bf 252} 115--139
% 
\bibitem{Kippenhahn1951} R.~Kippenhahn (1951)
{\em \"Uber den Wertevorrat einer Matrix},
Math Nachr {\bf 6} 193--228
%
\bibitem{KrupnikSpitkovsky2006} N.~Krupnik and I.\,M.~Spitkovsky (2006)
{\em Sets of matrices with given joint numerical range},
Linear Algebra Appl {\bf 419} 569--585
%
\bibitem{LiPoon2000} C.-K.~Li and Y.-T.~Poon (2000)
{\em Convexity of the joint numerical range},
SIAM J Matrix Anal A {\bf 21} 668--678
%
\bibitem{LoewyTam1986} R.~Loewy and B.-S.~Tam (1986)
{\em Complementation in the face lattice of a proper cone},
Linear Algebra Appl {\bf 79} 195--207
%
\bibitem{Magron-etal2015} V.~Magron, D.~Henrion, and J.-B.~Lasserre (2015)
{\em Semidefinite approximations of projections and polynomial images 
of semialgebraic sets}, 
SIAM J Optimiz {\bf 25} 2143--2164
%
\bibitem{Mueller2010} V.~M{\"u}ller (2010) 
{\em The joint essential numerical range, compact perturbations, and 
the Olsen problem}, Stud Math {\bf 197} 275--290
%
\bibitem{Netzer-etal2010} T.~Netzer, D.~Plaumann, 
and M.~Schweighofer (2010)
{\em Exposed faces of semidefinitely representable sets},
SIAM J Optimiz {\bf 20} 1944--1955
%
\bibitem{Ocko-etal2011}
S.\,A.~Ocko, X.~Chen, B.~Zeng, B.~Yoshida, Z.~Ji, M.\,B.~Ruskai, 
and I.\,L.~Chuang (2011)
{\em Quantum codes give counterexamples to the unique preimage 
conjecture of the N-representability problem},
Phys Rev Lett {\bf 106} 110501
%
\bibitem{Paulsen2002} V.\,I.~Paulsen (2002)
{\em Completely Bounded Maps and Operator Algebras},
New York: Cambridge University Press
%
\bibitem{RamanaGoldman1995} M.~Ramana and A.\,J.~Goldman (1995)
{\em Some geometric results in semidefinite programming},
J Global Optim {\bf 7} 33--50
%
\bibitem{Rockafellar1970} R.\,T.~Rockafellar (1970)
{\em Convex Analysis}, Princeton: Princeton University Press
%
\bibitem{Rodman-etal2016} L.~Rodman, I.\,M.~Spitkovsky, A.~Szko{\l}a, 
and S.~Weis (2016)
{\em Continuity of the maximum-entropy inference: Convex geometry and 
numerical ranges approach},
J Math Phys {\bf 57} 015204
%
\bibitem{RostalskiSturmfels2013} P.~Rostalski and B.~Sturmfels (2013)
{\em Dualities}, 203--249 of G.~Blekherman, P.~Parrilo, and R.~Thomas,
Eds., {\em Semidefinite Optimization and Convex Algebraic Geometry},
MOS-SIAM series on optimization,
Philadelphia: SIAM 
%
\bibitem{Sanyal-etal2011} R.~Sanyal, F.~Sottile, and B.~Sturmfels (2011)
{\em Orbitopes}, Mathematika {\bf 57} 275--314
%
\bibitem{Schneider2014} R.~Schneider (2014) 
{\em Convex bodies: The Brunn-Minkowski theory}, 2nd ed., 
New York: Cambridge University Press
%
\bibitem{SinnSturmfels2015} R.~Sinn and B. Sturmfels (2015)
{\em Generic spectrahedral shadows}, 
SIAM J Optimiz {\bf 25} 1209--1220
%
\bibitem{Szymanski-etal2016} K.~Szyma\'nski, S.~Weis, and K.~\.Zyczkowski 
(submitted)
{\em Classification of joint numerical ranges of three hermitian matrices 
of size three}, {\tt arXiv:1603.06569 [math.FA]}
%
\bibitem{Tam1985} B.-S.~Tam (1985)
{\em On the duality operator of a convex cone},
Lin Alg Appl {\bf 64} 33--56
%
\bibitem{Weis2011} S.~Weis (2011)
{\em Quantum convex support},
Linear Algebra Appl {\bf 435} 3168--3188
%
\bibitem{Weis2012a} S.~Weis (2012)
{\em A note on touching cones and faces},
J Convex Anal {\bf 19} 323--353
%
\bibitem{Weis2012b} S.~Weis (2012)
{\em Duality of non-exposed faces},
J Convex Anal {\bf 19} 815--835
%
\bibitem{WeisKnauf2012} S.~Weis and A.~Knauf (2012)
{\em Entropy distance: New quantum phenomena},
J Math Phys {\bf 53} 102206
%
\bibitem{Weis2014} S.~Weis (2014)
{\em Continuity of the maximum-entropy inference},
Commun Math Phys {\bf 330} 1263--1292
%
\bibitem{Ziegler1995} G.\,M.~Ziegler (1995)
{\em Lectures on Polytopes}, New York: Springer-Verlag
%
\end{thebibliography}

%
%
%
%%%%%%%%%%%%%%%%%%%%%%%%%%%%%%%%%%%%%%%%%%%%%%%%%%%%%%%%%%%%%%%%%%%%%%%%%%%%
%%%%%%%%%%%%%%%%%%%%%%%%%%%%%%%%%%%%%%%%%%%%%%%%%%%%%%%%%%%%%%%%%%%%%%%%%%%%
%%%%%%%%%%%%%%%%%%%%%%%%%%%%%%%%%%%%%%%%%%%%%%%%%%%%%%%%%%%%%%%%%%%%%%%%%%%%
%%%%%%%%%%%%%%%%%%%%%%%%%%%%%%%%%%%%%%%%%%%%%%%%%%%%%%%%%%%%%%%%%%%%%%%%%%%%
%%%%%%%%%%%%%%%%%%%%%%%%%%%%%%%%%%%%%%%%%%%%%%%%%%%%%%%%%%%%%%%%%%%%%%%%%%%%
%
\vspace{1cm}
\parbox{12cm}{%
Stephan Weis\\
e-mail: {\tt maths@stephan-weis.info}\\[.5\baselineskip]
Departamento de Matemática\\
Instituto de Matemática, Estatística e Computação Científica\\
Universidade Estadual de Campinas\\
Rua Sérgio Buarque de Holanda, 651\\
Campinas-SP, CEP 13083-859\\
Brazil\\[.5\baselineskip]
Centre for Quantum Information and Communication\\
Ecole Polytechnique de Bruxelles\\
Université Libre de Bruxelles\\
50 av. F.D. Roosevelt - CP165/59\\
B-1050 Bruxelles\\ 
Belgium}
%
%%%%%%%%%%%%%%%%%%%%%%%%%%%%%%%%%%%%%%%%%%%%%%%%%%%%%%%%%%%%%%%%%%%%%%%%%%%%
%%%%%%%%%%%%%%%%%%%%%%%%%%%%%%%%%%%%%%%%%%%%%%%%%%%%%%%%%%%%%%%%%%%%%%%%%%%%
%%%%%%%%%%%%%%%%%%%%%%%%%%%%%%%%%%%%%%%%%%%%%%%%%%%%%%%%%%%%%%%%%%%%%%%%%%%%
%%%%%%%%%%%%%%%%%%%%%%%%%%%%%%%%%%%%%%%%%%%%%%%%%%%%%%%%%%%%%%%%%%%%%%%%%%%%
%%%%%%%%%%%%%%%%%%%%%%%%%%%%%%%%%%%%%%%%%%%%%%%%%%%%%%%%%%%%%%%%%%%%%%%%%%%%
%
\end{document}